\documentclass[12pt]{elsarticle}

\usepackage{latexsym,amsfonts,amsmath,theorem,amssymb}
\usepackage{graphicx}

\def\be{\begin{equation}}
\def\bea{\begin{eqnarray*}}
\def\ee{\end{equation}}
\def\eea{\end{eqnarray*}}
\def\ba{\begin{array}}
\def\ea{\end{array}}
\def\bi{\begin{itemize}}
\def\ei{\end{itemize}}

\newtheorem{lem}{Lemma}
\newtheorem{pro}{Proposition}
\newtheorem{cor}{Corollary}

\def\ZZ{{{\rm Z}\kern-.4em{\rm Z}}}
\def\RR{{{\rm I}\kern-.2em{\rm R}}}
\def\NN{{{\rm I}\kern-.2em{\rm N}}}
\def\TT{{{\rm T}\kern-.5em{\rm T}}}
\def\CC{{{\rm I}\kern-.5em{\rm C}}}

\def\mX{\mathcal{X}}
\def\tmX{\tilde{\mathcal{X}}}
\def\mJ{\mathcal{J}}

\def\mP{\mathcal{P}}

\def\hbu{\hat{\mathbf{u}}}
\def\hbw{\hat{\mathbf{w}}}
\def\bfj{\mathbf{j}}
\def\bfi{\mathbf{i}}
\def\bfa{\mathbf{a}}
\def\bfb{\mathbf{b}}
\def\bfc{\mathbf{c}}
\def\bfs{\mathbf{s}}
\def\bfz{\mathbf{z}}
\def\bfW{\mathbf{W}}
\def\bfV{\mathbf{V}}
\def\bfU{\mathbf{U}}
\def\tf{\tilde{f}}
\def\tU{\tilde{U}}
\def\tW{\tilde{W}}
\def\tX{\tilde{X}}
\def\tK{\tilde{K}}

\def\tx{\tilde{x}}
\def\tu{\tilde{{u}}}
\def\tw{\tilde{\mathrm{w}}}
\def\tga{\tilde{\gamma}}
\def\tph{\tilde{\phi}}

\def\hu{\hat{\mathrm{u}}}
\def\hH{\hat{H}}

\def\homega{{\hat{\omega}}}
\def\ww{\mathrm{w}}
\def\vv{\mathrm{v}}
\def\uu{\mathrm{u}}
\journal{XXX}

\begin{document}
\begin{frontmatter}

\title{Hilbert function space splittings on domains with infinitely many variables}
\author[bonn,scai]{Michael Griebel}
\ead{griebel@ins.uni-bonn.de}
\author[bonn]{Peter Oswald}
\ead{agp.oswald@gmail.com}

\address[bonn]{Institute for Numerical Simulation (INS),
University of Bonn,
Wegelerstr. 6-8,
D-53211 Bonn, Germany}
\address[scai]{Fraunhofer-Institut f\"ur Algorithmen und Wissenschaftliches Rechnen SCAI, Schloss Birlinghoven, D-53754 Sankt
Augustin, Germany}
\date{}
\begin{abstract}
We present an approach to defining Hilbert spaces of functions depending on infinitely many  
variables or parameters, with emphasis on a weighted tensor product construction based on stable space splittings. The construction has been used in an exemplary way for guiding dimension- and scale-adaptive algorithms in application areas such as statistical learning theory, reduced order modeling, and information-based complexity. We prove results on compact embeddings, norm equivalences, and the estimation of $\epsilon$-dimensions. A new condition for the equivalence of
weighted ANOVA and anchored norms is also given. 
\end{abstract}
\end{frontmatter}

\section{Introduction}\label{sec0}
Functions depending on infinitely many real variables have been studied in different fields such as, e.g., stochastic processes, measure theory, stochastic 
and parametric PDEs, uncertainty quantification, and information-based complexity theory. A rich source of such functions is provided by nonlinear functionals
$F:\,K\subset U\to \mathbb{R}$ defined on a subset of a separable Banach space $U$. Such $F$ can often be parametrized by representing $K$ in the form
$ K=\{\sum_{k\in \mathbb{N}} x_ku_k:\; x=(x_1,x_2,\ldots)\in \mathcal{X}\subset \mathbb{R}^\infty\} $, with $\{u_k\}$ an appropriately chosen generating system of elements in $U$.
Then, $f(x):=F(\sum_{k\in \mathbb{N}} x_ku_k)$ is studied as function of infinitely many variables over the domain $\mathcal{X}\subset \mathbb{R}^\infty$ instead of $K\subset U$.

To classify functions of infinitely many real variables and to quantify their properties, various function spaces have been introduced following the traditional 
constructions for functions of finitely many variables, see the examples in Section \ref{sec1}. Often they arise via a tensor product construction from function spaces over one- or low-dimensional domains represented as direct sum of their subspaces. On the one hand, the use of tensor product techniques is bound to implicit assumptions on the separability or anisotropy of the functions under consideration, and to an ordering of the variables by significance. On the other hand, tensor product structures are instrumental when it comes to the efficiency of computational algorithms for high-dimensional problems. In some cases, spaces of functions depending on infinitely many variables arise naturally, e.g.,
in uncertainty quantification when input random fields are parametrized, or in control problems when the control belongs to an infinite-dimensional space. 
Alternatively, they are often constructed as suitable limits of their counterparts on domains $X^d\subset \mathbb{R}^d$ naturally embedded into $\mathcal{X}$, and serve as oracle spaces to analyze
and control the complexity of approximation algorithms in function spaces over $X^d$ as $d\to \infty$. In many case studies, the focus is on beating the curse of dimension by proving results that hold with estimates uniform in $d$, and the key to success is to establish a related stability result for a space on $\mathcal{X}\subset\mathbb{R}^\infty$. 

In this paper we consider separable Hilbert function spaces described in terms of decomposition norms. Our main goal is to establish a simple, yet flexible Hilbert space framework for studying functions on $\mathcal{X}$ which potentially implies robust results 
for the associated classes of spaces obtained by restriction to $X^d$, and covers most of the existing approaches. 
We will focus our attention on a construction based on Hilbert function spaces $H_k$ whose elements are functions $f_k:\, X_k\to \mathbb{R}$ of one variable
$x_k\in X_k\subset \mathbb{R}$, $k=1,2,\ldots$. 
Each of these $H_k$ is
split into a direct sum of closed subspaces $W_{j,k}$, where $j$ belongs to a finite or infinite index set $\mJ_k\subset \mathbb{Z}_+$ with $0\in \mJ_k$. For reasons of simplicity, we assume that $W_{0,k}=\mathrm{span}\{1\}$ is the space of all constant functions on $X_k$, with the normalization $\|1\|_{H_k}=1$. 
If $H_k=L_2(X_k,\mu_k)$, where $\mu_k$ is some probability measure over $X_k$, then direct sum splittings into one-dimensional $W_{j,k}$ arise from complete orthonormal systems (CONS) or Riesz bases. Splittings into orthogonal subspaces 
$$
L_2(X_k,\mu_k)=\bigoplus_{j\in\mJ_k} W_{j,k}, \qquad \|f_k\|_{L_2}^2 = \sum_{j\in\mJ_k} \|\ww_{j,k}\|_{L_2}^2, \qquad f_k = \sum_{j\in\mJ_k} \ww_{j,k},\quad \ww_{j,k}\in W_{j,k},
$$
and extensions to the case of non-direct, redundant space decompositions, related to so-called stable space splittings \cite{PO} or fusion frames \cite{CK}, a generalization of frames, are of interest as well.

Consistent with the product measure construction, we now introduce a whole family of Hilbert spaces of functions of infinitely many variables via
weighted Hilbert sums of tensor products of the $W_{j,k}$. 
Let  $\mJ$ denote the set of all index sequences $\bfj=(j_1,j_2,\ldots)\in \mJ_1\times \mJ_2\times \ldots$
with finite support $\omega_{\bfj}=\mathrm{supp}(\bfj) :=\{k\in\mathbb{N}:\;j_k>0\}$. The 
family $\bfW$ consists of tensor product Hilbert spaces 
$$
W_{\bfj}=\bigotimes_{k\in \omega_{\bfj}} W_{j_k,k},\qquad \bfj\in\mJ,
$$
whose elements can be interpreted as functions of finitely many variables (namely, the variables $x_k\in X_k$ with index $k$ in the support $\omega_{\bfj}$ of $\bfj$),
with the usual tensor product norm induced by
$$
\|\otimes_{k\in \omega_{\bfj}} \ww_{j_k,k}\|_{\omega_{\bfj}} := \prod_{k\in \omega_{\bfj}} \|\ww_{j_k,k}\|_{H_k}.
$$
Then $H^\infty_{\bfW,\bfa}$ is the set of all functions $f$ of infinitely many variables
admitting a unique decomposition $f = \sum_{\bfj\in \mJ} \ww_{\bfj}$ with $\ww_{\bfj}\in W_{\bfj}$, and finite norm defined by
\be\label{Norm}
\|f\|^2_{\bfW,\bfa} := \sum_{\bfj\in\mJ} a_{\bfj} \|\ww_{\bfj}\|_{\omega_{\bfj}}^2 < \infty.
\ee
The non-negative weight sequence $\bfa=\{a_{\bfj}\}_{\bfj\in\mJ}$ allows us to penalize the importance or to suppress the influence of certain components $\ww_{\bfj}$ or coordinate directions by letting $a_{\bfj}\to\infty$ or by setting $a_{\bfj}=0$, respectively. A more detailed definition will be given in the next section. The precise meaning of convergence of the series representation $f = \sum_{\bfj\in \mJ} \ww_{\bfj}$ and the natural domain of definition $\mathcal{X}$ of the functions $f\in H^\infty_{\bfW,\bfa}$ is not our main focus here, it will be elaborated on merely in connection with some examples. Since $H^\infty_{\bfW,\bfa}$ is the closure of all $f$ represented as the sum of finitely many $\ww_{\bfj}$ under the norm (\ref{Norm}), most of the properties of $H^\infty_{\bfW,\bfa}$ can be studied for
finitely representable $f$ for which convergence is not an issue, followed by a standard limit argument. 

The contributions of this paper are as follows. We first give a more detailed introduction to the spaces $H^\infty_{\bfW,\bfa}$ and their restrictions
$H^d_{\bfW,\bfa}$ to $X^d=X_1\times \ldots \times X_d$, given as the closed subspace of $H^\infty_{\bfW,\bfa}$ consisting of all 
$f$ of the form $f= \sum_{\bfj\in \mJ^d} \ww_{\bfj}$, where $\mJ^d=\mJ_1\times\ldots \times \mJ_d$. In other words, these are Hilbert decomposition spaces of functions depending on the first $d$ variables $x_k$, $k=1,\ldots,d$. This is done in Section \ref{Def} while in Section \ref{Ex} we review  the use of such spaces in previous research.
In Section \ref{sec2}, we give a few general statements on compact embeddings, $\epsilon$-dimensions, and the 
relationship between the spaces $H^\infty_{\bfW,\bfa}$ and $H^d_{\bfW,\bfa}$ as $d\to \infty$. We also study a particular case of redundant
decompositions, namely when the $H_k$ are equipped with an increasing ladder of subspaces $V_{j,k}=W_{0,k}+\ldots+W_{j,k}$, $j\in \mJ_k$, and 
introduce spaces  $H^\infty_{\bfV,\bfa}$  related to representations $f =   \sum_{\bfj\in \mJ} \vv_{\bfj}$, $\vv_{\bfj}\in V_{\bfj}$, with respect to a similarly constructed family $\bfV = \{V_{\bfj}:=\otimes_{k\in \omega_{\bfj}} V_{j_k,k}\}_{\bfj\in\mJ}$. Moreover, we compare the spaces $H^\infty_{\bfW,\bfa}$ and $H^\infty_{\bfV,\bfa}$.
Section \ref{sec3} shows the use of these families of spaces in concrete situations. In particular, we give estimates for $\epsilon$-dimensions
if the subspace family $\bfW$ and the weight sequence $\bfa$ are such that $H^\infty_{\bfW,\bfa}$ mimics a weighted tensor product of Sobolev
spaces $H^{s_k}(X_k)$. 
Next we show improved conditions on
the equivalence of norms associated with ANOVA and anchored decompositions that attracted some attention in connection with high-dimensional integration.
Finally, we consider the approximation problem for maps $F:\, K\subset U \to \tK \subset \tU$ between infinite-dimensional spaces in form of the closely related problem
of approximating functions $f:\,\mathcal{X} \to \tilde{\mathcal{X}}$ between parameter domains of infinitely many variables in the language of $H^\infty_{\bfW,\bfa}$ spaces.

There definitely is also need of studying spaces of functions of infinitely many variables that do not fit the Hilbert space setting considered in this paper. We refer to \cite{CD,D,ScGi} where this is highlighted in connection with optimal approximation procedures for stochastic PDEs, parameter-dependent problems, and reduced order modeling.
The extension to approximation spaces in an $L_p$ setting is a natural possibility, as is the incorporation of tensor products of Banach spaces.

\section{Definitions and Examples} \label{sec1}
\subsection{Definitions and Notation}\label{Def}
We first give the detailed definition of spaces $H^\infty_{\bfU,\bfa}$ for rather general subspace families $\bfU$ that contains the definition of the spaces $H^\infty_{\bfW,\bfa}$ and 
$H^\infty_{\bfV,\bfa}$ mentioned in Section \ref{sec0} as particular cases. 
\smallskip 
Our departure point is a sequence of Hilbert spaces $H_k$, $k=1,2,\ldots$, each of which comes with a collection of its closed subspaces
$\{U_{j,k}\}_{j\in\mJ_k}$ whose linear span is dense in $H_k$. The index sets $\mJ_k$ are non-empty subsets of $\mathbb{Z}_+$. Without loss of generality,
we may assume that $\mJ_k=\{0,1,\ldots,K_k\}$ for some finite $K_k\in \mathbb{Z}_+$ or $\mJ_k=\mathbb{Z}_+$ (then $K_k=\infty$). We assume that
elements of $H_k$ can be interpreted as functions $f_k:\,X_k\to \mathbb{R}$ of one variable $x_k$ with fixed domains $X_k\subset \mathbb{R}$. The assumption that the $H_k$ are functions of one real variable is not essential but we keep it for convenience. 

Concerning the subspaces $U_{j,k}$, we postulate that always $U_{0,k}=\mathrm{span}\{1\}$ 
with normalization $\|1\|_{H_k}=1$. 
To prepare for the upcoming discussion, we introduce the Hilbert sum 
$$
\hH_k=\sum_{j\in \mJ_k} U_{j,k}:= \{\hu_k:=(\uu_{0,k},\uu_{1,k},\ldots):\;\; \uu_{j,k}\in U_{j,k},\;j\in \mJ_k;\;\; \|\hu_k\|_{\hH_k} <\infty\},
$$
with the decomposition norm defined by
$$
\|\hu_k\|^2_{\hH_k}:=\sum_{j\in\mJ_k}\|\uu_{j,k}\|_{H_k}^2.
$$
The "sequence" space $\hat{H}_k$ can be related to $H_k$ by the following construction: For any $f_k\in \mathrm{span}\{U_{j,k}\}_{j\in\mJ_k}$, we define
$$
|f_k|^2_k := \inf_{f_k=\sum_j \uu_{j,k}} \|\hu_k\|_{\hH_k}^2,
$$
where the infimum is taken with respect to all $\hu_k\in\hH_k$ with finitely many $\uu_{j,k}\neq 0$. It is not hard to show that $|\cdot|_k$ is a 
Hilbert semi-norm on $\mathrm{span}\{U_{j,k}\}_{j\in\mJ_k}$, and that 
$$
|f_k|_k \le \|f_k\|_{H_k},\qquad f_k\in \mathrm{span}\{U_{j,k}\}_{j\in\mJ_k}.
$$
If a complementing lower bound $c_k\|f_k\|^2_{H_k}\le |f_k|^2_k$ holds with a constant $c_k>0$ independently of $f_k\in \mathrm{span}\{U_{j,k}\}_{j\in\mJ_k}$, then
$|\cdot|_k$ is a Hilbert norm, and the closure of $\mathrm{span}\{U_{j,k}\}_{j\in\mJ_k}$ under this norm leads to the same $H_k$, with spectrally
equivalent scalar products. This is the situation covered by the concept of stable space splittings \cite{PO} or the more recently introduced fusion frames
\cite{CK}. Both approaches contain the case of weighted Hilbert sums which is relevant for us: Given any sequence $a_k=\{a_{j,k}\ge 0\}_{j\in\mJ_k}$, we define
the space $\hH_{k,a_k}$ by the weighted sequence norm
$$
\|\hu_k\|^2_{\hH_k,a_k}:=\sum_{j\in\mJ_k} a_{j,k}\|\uu_{j,k}\|_{H_k}^2.
$$
It is silently assumed that a zero weight $a_{j,k}=0$ for some $j$ implies that the space $\hH_{k,a_k}$ consists only of sequences with $\uu_{j,k}=0$. 
As before, one can introduce a Hilbert semi-norm $|\cdot|_{k,a_k}$ which is now not necessarily related to the original $\|\cdot\|_{H_k}$. It turns out that the desirable norm property for $|\cdot|_{k,a_k}$ may hold for some $a_k$ but not for others. If it holds for some family of
weight sequences $a_k=\{a_{j,k}\}$ then the closure process with norm $|\cdot|_{k,a_k}$ leads to
a new scale of Hilbert spaces $H_{k,a_k}$ that, depending on the chosen subspace family $\{U_{j,k}\}$, may or may not have a natural interpretation. 
Note that in order to guarantee the norm property, one often needs a certain growth of $a_{j,k}$
as $j\to \infty$. Thus, due to the unconditionality of all summation processes, by reordering the $U_{j,k}$ we could silently assume monotonicity of the weight sequence for $j\ge 1$.

A simpler to analyze but important situation is when the subspaces $U_{j,k}$ form a Riesz system in $H_k$ in which case we denote them by
$W_{j,k}$. By definition, this means that any $f_k\in H_k$ has a unique series representation $f_k=\sum_{j\in\mJ_k} \ww_{j,k}$, $\ww_{j,k}\in W_{j,k}$, which converges unconditionally
in $H_k$ such that, with absolute constants $0<c_k\le C_k<\infty$, it holds
\be\label{RieszB}
c_k|f_k|_k^2 \le \|f_k\|_{H_k}^2 \le C_k|f_k|_k^2.
\ee
Here, $|\cdot|_k$ is defined with respect to the family $\{W_{j,k}\}_{j\in\mJ_k}$. Thus,  there is no need for the infimum in the norm definition
due to the uniqueness of decompositions. In this case, the original $H_k$ norm can be replaced by the spectrally equivalent norm $|\cdot|_k$ if needed. Moreover, if the Riesz system is orthogonal, i.e., if $H_k=\oplus_{j\in\mJ_k} W_{j,k}$, then these norms are identical. If in addition
all $W_{j,k}$ are one-dimensional, we can identify the above introduced $H_{k,a_k}$ with weighted $\ell_2$ spaces.


\smallskip
Next, we introduce tensor product spaces of functions depending on infinitely many variables $x=(x_1,x_2,\ldots)\in \mathcal{X}\subset X^\infty :=X_1\times X_2\times \ldots$. We again define  "sequence" spaces $\hH^\infty_{\bfU,\bfa}$ first, and then, by a similar summation and closure process as above, function spaces $H^\infty_{\bfU,\bfa}$. 
To this end, let $\mJ$ denote the set of all  index sequences 
$\bfj=(j_1,j_2,\ldots)$, $j_k\in\mJ_k$, with finite support denoted by
$$
\omega_{\bfj}=\mathrm{supp}(\bfj) :=\{k\in\mathbb{N}:\;j_k>0\}\in \mP_f(\mathbb{N}),\qquad |\omega_{\bfj}|=|\bfj|_0 <\infty.
$$  
Here, $\mP_f(\mathbb{N})$ denotes the set of all finite subsets of $\mathbb{N}$, $|\omega|$ is the number of elements in a set $\omega$,
and $|\bfj|_0$ is the number of non-zero entries of the vector $\bfj$. 

Then the family $\bfU$  consists  of the tensor product spaces $U_{\bfj}:=\bigotimes_{k\in\omega_{\bfj}} U_{j_k,k}$ whose elements are identified with functions on the tensor product domains 
$$
X^{\omega_{\bfj}}:=\times_{k\in\omega_{\bfj}} X_k, \qquad \bfj\in\mJ,
$$ 
depending only on the variables $x_k$ with $k\in \omega_{\bfj}$ via the formula
$$
\otimes_{k\in\omega_{\bfj}} \uu_{j_k,k} = \prod_{k\in\omega_{\bfj}} \uu_{j_k,k}(x_k),\qquad \uu_{j_k,k}\in U_{j_k,k},
$$
for elementary tensors. The standard tensor product norm on $U_{\bfj}$ is induced by setting 
$$
\|\otimes_{k\in\omega_{\bfj}} \uu_{j_k,k}\|_{\omega_{\bfj}} := \prod_{k\in\omega_{\bfj}} \|\uu_{j_k,k}\|_{H_k}.
$$
Indexing this norm by $\omega_{\bfj}$ is appropriate, as it depends only on the support of $\bfj$, and not on the specific entries $j_k\neq 0$.
Now we are ready to define the "sequence" spaces
\be\label{hHU}
\hH_{\bfU,\bfa}=\sum_{\bfj\in \mJ} a_{\bfj}U_{\bfj}:= \left\{\hbu=\{\uu_{\bfj}:\;\uu_{\bfj}\in U_{\bfj},\;\bfj\in \mJ\}:\;\; \|\hbu\|_{\bfU,\bfa} <\infty\right\},
\ee
where $\bfU:=\{U_{\bfj}:\,\bfj\in\mJ\}$ denotes the underlying family of tensor product spaces while the decomposition norm is defined by
$$
\|\hbu\|^2_{\bfU,\bfa}:=\sum_{\bfj\in\mJ} a_{\bfj}\|\uu_{\bfj}\|_{\omega_{\bfj}}^2.
$$
The weight sequence $\bfa$ is always non-negative. Moreover, $a_{\bfj}=0$ for some $\bfj\in \mJ$ implies that the corresponding component $\uu_{\bfj}$ in $\hbu$
vanishes or, equivalently, that the associated space $U_{\bfj}$ is dropped from $\bfU$. As a consequence, any summations with respect
to $\mJ$ are silently to be taken with respect to the support $\omega_{\bfa}:= \{\bfj\in \mJ:\;a_{\bfj}>0\}$ of $\bfa$. 

As before, the "sequence" space $\hH_{\bfU,\bfa}$ can be turned into a function space.  For any $f\in \mathrm{span} \bfU$, define a semi-norm by setting
\be\label{HUnorm}
|f|^2_{\bfU,\bfa} := \inf_{f=\sum_{\bfj} u_{\bfj}} \sum_{\bfj\in\mJ}  a_{\bfj}\|\uu_{\bfj}\|_{\omega_{\bfj}}^2,
\ee
where the infimum is taken with respect to all possible finite decompositions of $f\in \mathrm{span} \bfU$. Note that $\mathrm{span} \bfU$ consists of functions well-defined on $X^\infty$, with the property that they are non-constant only with respect to a finite number of variables $x_k$. 
We now make the assumption that (\ref{HUnorm}) defines a norm on $\mathrm{span}\bfU$, i.e., that $|f|_{\bfU,\bfa}=0$ for $f\in\mathrm{span}\bfU$ implies $f=0$. Then, by definition, $H^\infty_{\bfU,\bfa}$ is introduced as the closure of $\mathrm{span}\bfU$ with respect to this norm, and we are allowed to
use the symbol $\|\cdot\|_{\bfU,\bfa}$ for the extension of the norm $|\cdot|_{\bfU,\bfa}$ to this closure. 

The norm property is automatic for any weight sequence $\bfa$ if $\bfU=\bfW$, i.e., when the $\{W_{j,k}\}_{j\in\mJ_k}$ are Riesz systems in $H_k$ for all $k\in\mathbb{N}$, see (\ref{RieszB}). In this case, the spaces $H^\infty_{\bfW,\bfa}$ and $\hH_{\bfW,\bfa}$ can be identified via an isometry, as each $\hbu\in H^\infty_{\bfW,\bfa}$
generates a unique equivalence class of Cauchy sequences in $\mathrm{span}\bfU$. This can be seen by taking any sequence of finite "partial sums" with entries from $\hbu$. Vice versa, any equivalence class of Cauchy sequences from $\mathrm{span}\bfU$ belonging to $H^\infty_{\bfW,\bfa}$ defines a unique $\hbu\in \hH_{\bfW,\bfa}$ such that its partial sum sequences belong to the same equivalence class. 
From now on, we use the notation 
$$ 
f = \sum_{\bfj\in \mJ} \ww_{\bfj}, \qquad \|f\|_{\bfW,\bfa}^2=\sum_{\bfj\in\mJ} a_{\bfj} \|\ww_{\bfj}\|_{\omega_{\bfj}}^2,
\qquad \hbw=\{\ww_{\bfj}:\;\ww_{\bfj}\in W_{\bfj},\;\bfj\in \mJ\}, 
$$
to express the relationship between $f\in H^\infty_{\bfW,\bfa}$ and the unique $\hbw\in \hH_{\bfW,\bfa}$ associated with it. 
If $\bfU\neq \bfW$, i.e., if some (or all) of the $H_k$ are equipped with redundant space splittings,  the weights for which (\ref{HUnorm}) defines a norm can in some cases be completely characterized, see Section \ref{sec2}. 
Although mathematically convenient, the implicit definition of $H^\infty_{\bfU,\bfa}$ via completion (implicit because equivalence classes of Cauchy sequences
in $\mathrm{span}\bfU$ are not functions) can sometimes be replaced by an intrinsic description as space of functions
defined on some $\mathcal{X}\subset X^\infty$, see \cite{Ar,GMR} for the case of tensor product reproducing kernel Hilbert spaces. However, we neglect this issue for most part of this paper. 

\smallskip 
After having defined the scales of weighted decomposition Hilbert spaces $H^\infty_{\bfU,\bfa}$ of functions of infinitely many
variables, 
we can identify various subspaces of interest for high-dimensional approximation. 
If $\bfU=\bfW$, this is simply done by restriction of the admissible sequences $\hbw$ or the weight sequence $\bfa$.
For instance, spaces of functions on $X^d$ are defined by requiring $\ww_{\bfj}=0$ whenever $\omega_{\bfj}\not\subset \{1,\ldots,d\}$, or equivalently by replacing $\bfa$
by its restriction to $\mJ^d$, i.e., by setting
$$
H^d_{\bfW,\bfa}:=H^\infty_{\bfW,\bfa^d},\qquad a_{\bfj}^d:=\left\{\ba{ll} a_{\bfj},& \omega_{\bfj}\subset \{1,\ldots,d\},\\
0, &\omega_{\bfj}\not\subset \{1,\ldots,d\}.\ea\right.
$$
In this case, the summation process for defining $f$ is usually written as $d$-fold summation,
$$
f(x):= \sum_{j_1=0}^{K_1}\cdots \sum_{j_d=0}^{K_d} \ww_{j_1,\ldots,j_d}(x), \quad x=(x_1,\ldots,x_d),
$$
by leaving out all unnecessary variables $x_k$ with $k>d$ and components $\ww_{\bfj}$, where $\omega_{\bfj}\not\subset \{1,\ldots,d\}$.
Questions on the behavior of approximation processes for functions from $H^d_{\bfW,\bfa}$ for large $d$, and asymptotically for $d\to \infty$
are intimately related to the  "limit" space $H^\infty_{\bfW,\bfa}$.

Some applications focus on approximation by functions of fewer variables, say, at most $m$-variables ($m\ge 1$, excluding the trivial case $m=0$ of approximation by global constants). In this case, one suppresses components $\ww_{\bfj}$ for which $|\omega_{\bfj}|>m$, and considers representations of the form
$$
f(x) = \ww_{\mathbf{0}} + \sum_{k\in\mathbb{N}} f_k(x_k) + \ldots + \sum_{(k_1,\ldots,k_m)\in \mathbb{N}^m} f_{k_1,\ldots,k_m}(x_{k_1},\ldots,x_{k_m}),
$$
where the functions $f_{k_1,\ldots,k_l}(x_{k_1},\ldots,x_{k_l})$ depend only on the indicated variables, and belong to the closure of the span of all those $W_{\bfj}$ for which $\omega_{\bfj}=\{k_1,\ldots,k_l\}$. 
To cover these applications, we would set all $a_{\bfj}$ with $|\omega_{\bfj}|>m$ to zero in the restriction process. A slightly more compact expression results if we index by subsets of $\mathbb{N}$:
\be\label{mVar}
f(x) = \sum_{l=0}^m \sum_{\omega\subset\mathbb{N}:\,|\omega|=l} f_\omega(x_\omega),\qquad f_\omega(x_\omega)=\sum_{\bfj \in\mJ:\, \omega_{\bfj}=\omega} 
\ww_{\bfj}(x_\omega),\quad x_\omega=x|_\omega\in X^\omega.
\ee
 
There is yet another option of defining a natural order in the summation processes to represent elements of $H^\infty_{\bfW,\bfa}$, namely by increasing
values of $a_{\bfj}$. Then one would define "partial sums"
$$
f_M:= \ww_{\mathbf{0}}+ \sum_{\mathbf{0}\neq\bfj\in\mJ:\,0<a_{\bfj}\le M} \ww_{\bfj},
$$ 
assuming that the respective index sets are finite, and considers $f$ as limit of $f_M$ for $M\to \infty$. The rationale of such an ordering is that large values of $a_{\bfj}$ mean penalization of the corresponding terms $\ww_{\bfj}$ in the representation
of $f$. Thus, such $\ww_{\bfj}$ can be considered small in size, and can be added later. In many of the examples below, the weights $a_{\bfj}$, $\bfj\neq\mathbf{0}$, are increasing with respect to the partial order on $\mJ$, i.e.,
$$
\mathbf{0}<\bfj\le \bfi \quad \Longrightarrow \quad a_{\bfj}\le a_{\bfi}.
$$
This monotonicity reflects  the traditional ordering in the $k$- and $j$-scales, where the dependence on variables $x_k$ with large $k$ is less significant.
Larger $j_k$ corresponds to higher frequency content in direction $x_k$ which is penalized by larger weights if smoothness in direction $x_k$ is assumed. This also leads to index sets for the summation of $f_M$ that are monotone with respect to the ordering in $\mJ$. 

Further examples, and classes of weight sequences that appear in certain applications, are given below. We finish with the remark that the definition of subspaces 
of $H^\infty_{\bfU,\bfa}$ for redundant families $\bfU\neq \bfW$ is a much more subtle issue related to the fact that subsets of frames may not be frames at all,
i.e., they may not inherit any nice properties of $H^\infty_{\bfU,\bfa}$. For the case $\bfU=\bfV$, see Sections \ref{sec2} and \ref{sec4}.

\subsection{Examples}\label{Ex}
To show that the definitions given in the previous section are meaningful, we compose a short list of examples from the literature, and interpret them within the family of $H^{\infty}_{\bfU,\bfa}$ spaces.
\subsubsection{Spaces related to Gaussian measures}\label{Ex1} Functions whose domain is a subset of a infinite-dimensional linear space are a natural source for our setting. We start with a classical, yet instructive example of function spaces on a separable Hilbert space $H$ related to the semi-group theory of
stochastic processes, see \cite{dP} for the basic theory, and \cite{Bo} for a more in-depth exposition.  
Given a symmetric positive definite trace class operator $Q:\, H\to H$, we can identify $H$ with $\ell^2$ using the eigensystem $\{e_k,\lambda_k\}$ associated with $Q$. Thus,  functions $f:\,H\to \mathbb{R}$ with argument $u=\sum_{k=1}^\infty x_k e_k\in H$ can be identified with functions $\tf:\, \mX=\ell^2\subset \mathbb{R}^\infty$ given by
$$
\tf(x)=f(\sum_{k=1}^\infty x_k e_k),\qquad x=(x_1,x_2,...)\in\ell^2,\quad\quad \|x\|_{\ell^2}=\|u\|_H^2.
$$
The centered Gaussian measure $\mu_{Q}$ with covariance operator 
$Q$ is then defined as the unique probability measure on $(H,\mathcal{B}(H))$ induced by the product Gauss measure
$\tilde{\mu}_{Q}(dx) = \otimes_{k=1}^\infty N_{0,\lambda_k}(dx_k)$
of univariate centered normal distributions with variance $\lambda_k$ on $\ell^2$.

The space $L_2(H,\mu_Q)$ is introduced via $L_2(\ell^2,\tilde{\mu}_Q)$ as usual. Then 
the tensor product Hermite polynomials
$$
h_{\bfj}(x)=\prod_{k\in \omega_{\bfj}} h_{j_k}(\frac{x_k}{\sqrt{\lambda_k}}),\qquad \bfj\in\mJ,
$$
can be shown to be a CONS in $L_2(\ell^2,\tilde{\mu}_Q)$, where $h_j(s)$ are the univariate Hermite polynomials. 

This fits our abstract scheme if the Hilbert spaces $H_k$ are identified with $L_2(\mathbb{R},N_{0,\lambda_k})$ with domain $X_k=\mathbb{R}$, and are split into
one-dimensional $W_{j,k}$ spanned by the Hermite polynomials $h_j(x_k/\sqrt{\lambda_k})$. With these ingredients at hand, by orthonormality,
we automatically have $L_2(\ell_2,\tilde{\mu}_Q)=H_{\bfW,\mathbf{1}}^\infty$,
where $\mathbf{1}$ denotes the weight sequence of all ones. Note here that ordering the decomposition of $\tf\in L_2(\ell_2,\tilde{\mu}_Q)$ by total polynomial degree, i.e.,
$$
\tf =\sum_{n=0}^\infty I_n(\tf),\qquad I_n(\tf)=\sum_{\bfj\in\mJ:\,|\bfj|_1=n} \ww_{\bfj}, \quad \ww_{\bfj}=(\tf,h_{\bfj})_{L_2}h_{\bfj}\in W_{\bfj},\quad \bfj\in\mJ,
$$
is called Wiener-Ito or Hermite polynomial chaos expansion, and represents yet another possibility of interpreting the summation process underlying the $H^\infty_{\bfW,\bfa}$ spaces. In this setting, $L_2$ based Sobolev spaces of positive smoothness
\bea
H^s(\ell^2,\tilde{\mu}_Q) &:= &\{\tf \in L_2(\ell_2,\tilde{\mu}_Q):\; \|\tf\|^2_{H^s}:=\sum_{n=0}^\infty (n+1)^{2s}\|I_n(\tf)\|_{L_2}^2<\infty\},\qquad s>0,
\eea
can be identified with $H_{\bfW,\bfa}^\infty$ spaces if the weight sequence $\bfa$ is defined by
$$
a_{\bfj}:= (1+ \sum_{k\in \omega_{\bfj}} \frac{j_k}{\lambda_k})^{2s},\qquad \bfj\in\mJ.
$$
We refer to \cite[Chapter 5]{Bo} for various equivalent definitions of $H^s(\ell^2,\tilde{\mu}_Q)$, $s\in \mathbb{N}$, and their extensions to $L_p$ based Sobolev spaces of Hilbert and Banach space valued functions over locally convex spaces equipped with a Gaussian measure.  

\subsubsection{Example 2: Decomposition spaces for uncertainty quantification and parameter identification}\label{Ex2} In the application fields of uncertainty quantification and parameter identification,
solutions to problems depending on a basic variable $\tilde{x} \in \tmX$ and infinitely many random or deterministic parameters $x=(x_1,x_2,\ldots)\in \mX$ are often modeled by expansions of the form
\be\label{ParaF}
f(\tilde{x};x)= \sum_{\bfj\in \mJ} f_{\bfj}(\tilde{x}) \phi_{\bfj}(x_{\omega_{\bfj}})
\ee
where $\phi_{\bfj}(x_{\omega_{\bfj}})=\prod_{k\in \omega_{\bfj}} \phi^{k}_{j_k}(x_k)$ are tensor product basis functions adapted to the problem at hand. In the setting of 
random parameter vectors $x$, assuming mutual independence of the random variables $x_k$ with underlying marginal distributions $\mu_k$ and uniformly bounded second moments, the univariate systems $\{\phi^{k}\}$ are chosen such that
they form complete orthogonal systems in $L_{2,\mu_k}(X_k)$, where $X_k\subset \mathbb{R}$ is identified with the support of $\mu_k$.
Often systems of orthogonal polynomials on intervals or discrete sets are adopted, see \cite[Section 10]{Sm},
\cite{KSXLST,WK,XK}, the Hermite polynomial expansion in Section \ref{Ex1} being a particular example. Similar decompositions are used in the case of deterministic parameter vectors \cite{CD}, here the measures $\mu_k$ do not have a probabilistic interpretation. Needless to mention that, in addition to orthogonal expansions and generalized polynomial chaos, more general decompositions can be used in (\ref{ParaF}), see, e.g., \cite{BCM,BCDM,BNT,CD,ScGi}. 

If the functions $f_{\bfj}$ are elements of a Hilbert space $\tilde{H}$ of functions on $\tX$, then the decomposition (\ref{ParaF}), now viewed as a $\tilde{H}$-valued function on $\mX$, generates an orthogonal splitting of 
$L_{2,\mu}(X^\infty,\tilde{H})$, where $\mu$ is the product measure on $\mX\subset X^\infty$
induced by the measures $\mu_k$. Taking $H_k=L_{2,\mu_k}(X_k,\tilde{H})$, $W_{j,k}=\{\phi^k_j(x_k)\tilde{u}:\;\tilde{u}\in \tilde{H}\}$, and introducing a weight sequence $\bfa$, creates a huge family of  Hilbert spaces $H_{\bfW,\bfa}(X^\infty,\tilde{H})$ of $\tilde{H}$-valued functions that can be used to 
classify the convergence behavior of linear and nonlinear approximation processes based on (\ref{ParaF}).
When only a single quantity of interest is modeled then the dependence on $\tx$ can be neglected, and we are in the setting of spaces $H^\infty_{\bfW,\bfa}$
of Section \ref{Def}, with $\bfW$ induced by the system $\{\phi_{\bfj}\}_{\bfj\in\mJ}$.

\subsubsection{Example 3: Tensor product reproducing kernel Hilbert spaces}\label{Ex3} Decomposition spaces based on tensor products of reproducing kernel Hilbert spaces (RKHS) go back to \cite{Ar} and arise in statistical learning theory \cite{CZ,GKKW,MMRTS,SC}
and data analysis where they serve as hypothesis spaces. They are also used for studying tractability of function approximation and integration problems, especially if standard information is considered, see, e.g., \cite{KSWW,W,WW2} and Section \ref{sec5}. Roughly speaking, a Hilbert space $H$ of functions $f: \,\Omega \to \mathbb{R}$ is a RKHS characterized by a positive definite kernel $\kappa:\,\Omega\times \Omega\to \mathbb{R}$ such that function evaluation at any $x\in\Omega$ is a continuous functional on $H$ represented by the kernel, i.e., $f(x)=(f,\kappa(x,\cdot))_H$ for all $f\in H$ and $x\in\Omega$. Least-squares regression problems with RKHS induced penalty terms can conveniently be solved using finite linear combinations of the kernel, we refer to Section \ref{sec6} for more details. 

Often, $\Omega$ is a compact subset of $\mathbb{R}^d$ with large $d$, 
and $H$ arises from a tensor product construction
of RKHS $H_k$ decomposed into a direct sum of subspaces.
In \cite{ES}, weighted decomposition spaces $H_{\gamma}$ were introduced by partitioning $\Omega$ into pairwise disjoint
regions $\Omega_j$, each equipped with a different hypothesis RKHS $H_j$ of functions $f_j:\,\Omega_j \to \mathbb{R}$. Then, by definition, $f\in H_{\gamma}$ if
$f|_{\Omega_j}=f_J\in H_j$ for all $j$, and 
$$
\|f\|_{H_{\gamma}}^2 := \sum_j \gamma_j \|f_j\|_{H_j}^2 <\infty.
$$
The parameters $\gamma_j>0$ are used to incorporate the relative importance of
the hypotheses on the local behavior of $f$ expressed by the choice of $H_j$. When combined with tensorization, we arrive at our setting.

In \cite{Gu,Wa}, a method called SS-ANOVA (smoothing spline analysis of variance) is described, where a RKHS $H=H_0\bigoplus H'$ is orthogonally decomposed into  a low-dimensional subspace $H_0$ with kernel $\kappa_0$, and a RKHS $H'$ with kernel $\kappa'$.  Often, $H_0$ is further split into one-dimensional subspaces $W_j$, $j=0,\ldots,K-1$, and we set $W_K=H'$ to arrive at the situation discussed in our paper, see \cite[Chapter 2 and 4]{Gu} or \cite{Wa} for concrete examples.  When multivariate problems on product domains $X^d=X_1\times \ldots X_d$ with $X_k\subset \mathbb{R}$ are discussed, the common strategy
in SS-ANOVA schemes \cite{Gu,Wa,Wa1} is to use tensorization based on coordinate RKHS spaces $H_k= \bigoplus_{j\in\mJ_k} W_{j,k}$ constructed as explained before, in line with our $H^d_{\bfW,\bfa}$ construction. 
One specific issue in SS-ANOVA applications is the choice of the weights $a_{\bfj}$, many of which are set to zero. This is done by an a priori or a posteriori model selection step. Moreover, some of the weights are adjusted when solving the regression problem computationally to address the bias-variance problem, see \cite{CZ,Gu,Wa}. 

\subsubsection{Example 4: Spectral decompositions of smoothness spaces based on CONS}\label{Ex4} If each of the Hilbert spaces $H_k$ is equipped with a CONS $\{e_{j,k}\}_{j\in\mJ_k}$, one can create orthogonal
decompositions $H_k=\bigoplus_{j\in\mJ_k} W_{j,k}$ with one-dimensional subspaces $W_{j,k}$ spanned by the individual elements $e_{j,k}$ of these CONS, or, as it would be natural for wavelet type systems, subspaces $W_{j,k}$ spanned by all orthogonal wavelet or semi-orthogonal prewavelet functions of level $j$. 
E.g., in \cite{DG} best linear approximation has been investigated in spaces $H_{\bfW,\bfa}^\infty$ based on such CONS splittings with weight sequences $\bfa$ that model varying finite order, mixed and isotropic Sobolev smoothness or analytic behavior on an infinite product domain $X^\infty$. Here, periodic intervals $X_k=\mathbb{T}$, equipped with the usual Lebesgue measure and the trigonometric CONS for each component space, and non-periodic
situations such as $X_k=[-1,1]$ (Lebesgue measure, Legendre system as CONS) and $X_k=\mathbb{R}$ (Gaussian measure, Hermite polynomials as CONS) have been considered.

For the case $X_k=\mathbb{T}$ and the trigonometric CONS, we also refer to \cite{DTU} which contains a detailed survey of approximation results
for various function spaces on $X^d$ for finite $d$, asymptotics for $d\to \infty$, and a short paragraph about $d=\infty$, including the Hilbert space setting of the present paper. Under similar assumptions, and for rather general classes of weights, the recent papers \cite{KMU,KSU} contain precise asymptotics and pre-asymptotics for approximation numbers of the natural embedding operators of such  $H_{\bfW,\bfa}^d$ spaces.

\section{Theoretical Results}\label{sec2}
We will state some facts about embeddings and rates of best linear approximation for the spaces $H^\infty_{\bfW,\bfa}$,
and discuss a particular instance of spaces based on redundant decompositions, namely $H^\infty_{\bfV,\bfa}$.
This is a relatively easy task, as our definition of these spaces is essentially based on the "sequence" spaces $\hH_{\bfW,\bfa}$ 
which can be viewed as spaces $\ell^2_{\bfa}(\mJ)$ with entries from $W_{\bfj}$, for which embeddings and best approximation processes are well-understood.

We state the first, almost obvious result on embeddings for spaces $H^\infty_{\bfW,\bfa}$ based on non-redundant decompositions. Recall that $\omega_{\bfa}$
denotes the support of $\bfa$.
\begin{pro}\label{prop1} Let the $\{W_{j,k}\}_{j\in\mJ_k}$ be Riesz systems in $H_k$ for any $k=1,2,\ldots$, and consider arbitrary
weight sequences $\bfa$ and $\bfb$.\\
a) We have $H^\infty_{\bfW,\bfa}\subset H^\infty_{\bfW,\bfb}$ if and only if $\omega_{\bfa} \subset \omega_{\bfb}$,
and if there is a constant $C<\infty$ such that $b_{\bfj}\le C^2a_{\bfj}$ for all $\bfj\in \omega_{\bfa}$. Moreover, inclusion automatically implies continuous 
embedding, with norm of the embedding operator $\le C$.\\
b) Under the conditions in a), the embedding is compact if and only if for any $\bfj\in \omega_{\bfa}$ the associated $W_{\bfj}$ is finite-dimensional, and if for any $\epsilon>0$ the set 
$$
\mJ_{\bfc,\epsilon}:= \{\bfj\in\mJ:\; c_{\bfj}\ge \epsilon^2 \}=\{\bfj\in\mJ\cap \omega_{\bfa}:\; b_{\bfj}\ge \epsilon^2 a_{\bfj}\}
$$
is finite. Here, $\bfc$ is defined by $c_{\bfj}=b_{\bfj}/a_{\bfj}$ if $a_{\bfj}>0$ and
$c_{\bfj}=0$ if $a_{\bfj}=0$.
\end{pro} 

The proof of Part a) of Proposition \ref{prop1} is a direct consequence of the definition of the spaces $H^\infty_{\bfW,\bfa}$, its Part b) is part of the proof of the following Proposition \ref{prop2} related to best linear approximation in $H^\infty_{\bfW,\bfa}$. 

We are particularly interested in the so-called $\epsilon$-dimension 
$n_\epsilon(H^\infty_{\bfW,\bfa},H^\infty_{\bfW,\bfb})$ which, for given $\epsilon>0$, is defined as the smallest $n$ such that there exists a linear subspace 
$M\subset H^\infty_{\bfW,\bfb}$ of dimension $\dim M\le n$ with the property
\be\label{neps}
\inf_{g\in M} \|f-g\|_{\bfW,\bfb} \le \epsilon  \|f\|_{{\bfW,\bfa}},\qquad \forall\; f\in H^\infty_{\bfW,\bfa}.
\ee
The $\epsilon$-dimension is the inverse function to the Kolmogorov $n$-width of the unit ball of $H^\infty_{\bfW,\bfa}$ in $H^\infty_{\bfW,\bfb}$.
Finiteness of $n_\epsilon(H^\infty_{\bfW,\bfa},H^\infty_{\bfW,\bfb})$ for all $\epsilon>0$ is equivalent to the statement of Part b) in Proposition \ref{prop1}.


\begin{pro}\label{prop2}
Using the notation in Part b) of Proposition \ref{prop1}, we have
$$
n_\epsilon(H^\infty_{\bfW,\bfa},H^\infty_{\bfW,\bfb}) = \sum_{\bfj\in \mJ_{\bfc,\epsilon}} \dim W_{\bfj},
$$
where  $ \dim W_{\bfj}= \prod_{k\in \omega_{\bfj}} d_{j_k,k}$  for $\mathbf{0}\neq \bfj\in\mJ$, with $d_{j,k}:=\dim W_{j,k}$ denoting the dimensions of the coordinate subspaces, and $\dim W_{\mathbf{0}}=1$.
\end{pro}

{\bf Proof}. Obviously,  taking the specific linear subspace 
$M_{\bfc,\epsilon}=\sum_{\bfj\in\mJ_{\bfc,\epsilon}} b_{\bfj}W_{\bfj}$ of $H^\infty_{\bfW,\bfb}$ and any $f=\sum_{\bfj\in\mJ} \ww_{\bfj}\in
H^\infty_{\bfW,\bfa}$ results in the estimate
\bea
\mathrm{dist}_{H^\infty_{\bfW,\bfb}}(f,M_{\bfc,\epsilon})^2 &=& \sum_{\bfj\not\in \mJ_{\bfc,\epsilon}} b_{\bfj} \|\ww_{\bfj}\|^2_{\omega_{\bfj}}
\le (\sup_{\bfj\not\in \mJ_{\bfc,\epsilon}} c_{\bfj}) \sum_{\bfj\not\in \mJ_{\bfc,\epsilon}} a_{\bfj} \|\ww_{\bfj}\|^2_{\omega_{\bfj}}\\
&\le& \epsilon^2 \sum_{\bfj\in \mJ} a_{\bfj} \|\ww_{\bfj}\|^2_{\omega_{\bfj}} = \epsilon^2 \|f\|_{\bfW,\bfa}^2.
\eea
This gives the upper estimate 
$$
n_\epsilon(H^\infty_{\bfW,\bfa},H^\infty_{\bfW,\bfb})\le \dim M_{\bfc,\epsilon}=\sum_{\bfj\in \mJ_{\bfc,\epsilon}} \dim W_{\bfj}.
$$
If a finite-dimensional linear subspace $M$ of $H^\infty_{\bfW,\bfb}$ does not contain $M_{\bfc,\epsilon}$ then there must be
at least one $\bfj \in\mJ_{\bfc,\epsilon}$ and a non-zero element $\ww^\ast_{\bfj}\in W_{\bfj}$ that is orthogonal to $M$. 
For this $\ww^\ast_{\bfj}$, considered as element in $H^\infty_{\bfW,\bfa}$, we have
$$
\mathrm{dist}_{H^\infty_{\bfW,\bfb}}(\ww^\ast_{\bfj},M)^2=\|\ww^\ast_{\bfj}\|^2_{\bfW,\bfb}=b_{\bfj}\|\ww^\ast_{\bfj}\|_{\omega_{\bfj}}^2 =c_{\bfj}\|\ww^\ast_{\bfj}\|^2_{\bfW,\bfa} >\epsilon^2\|\ww^\ast_{\bfj}\|^2_{\bfW,\bfa}.
$$
Thus, $M$ does not provide the approximation quality needed in (\ref{neps}), only subspaces $M\supset M_{\bfc,\epsilon}$ can do. This shows the
equality $n_\epsilon(H^\infty_{\bfW,\bfa},H^\infty_{\bfW,\bfb}) = \dim M_{\bfc,\epsilon}$, and establishes  Proposition \ref{prop2}. 

\medskip
Proposition \ref{prop2} reduces the estimation problem for $\epsilon$-dimensions and related widths to a combinatorial optimization problem once $\bfa$, $\bfb$, and the dimensions $d_{j,k}$ are given.
In principle, this can be done computationally if good estimates for these sequences are known. 
If the $W_{j,k}$ are one-dimensional subspaces generated by complete orthonormal systems $\{e_{0,k}=1,e_{j,k}: \,j\ge 1\}$ in $H_k$, $k\in \mathbb{N}$,
then $d_{j,k}=1$ and 
$$
n_\epsilon(H^\infty_{\bfW,\bfa},H^\infty_{\bfW,\bfb})=|\mJ_{\bfc,\epsilon}|, 
$$
and we can use counting arguments or volume estimates if an explicit description of $\mJ_{\bfc,\epsilon}$ is available.
The compactness assumption says that $c_{\bfj}=b_{\bfj}/a_{\bfj}$ should tend to zero on $\omega_{\bfa}$. To achieve slow growth of $n_\epsilon(H^\infty_{\bfW,\bfa},H^\infty_{\bfW,\bfb})$ the decay of the $c_{\bfj}$ has to be the faster, the larger the dimensions $d_{\bfj}$ and the number of involved variables $|\bfj|_0$ become. One option to enforce this is to decrease the weights $b_{\bfj}$ for $\bfj$ with large support or with large $\max(\mathrm{supp}(\bfj))$ which de-emphasizes the approximation quality in higher dimensions. Alternatively, one can increase the weights $a_{\bfj}$ which can often be interpreted as requiring more smoothness for the function to be  approximated. 

\medskip
Next, we draw some conclusions for the $d$-dimensional counterparts of the $H_{\bfW,\bfa}(X^\infty)$ spaces
$$
H^d_{\bfW,\bfa}\equiv H_{\bfW,\bfa}(X^d) :=\sum_{\bfj\in \mJ^d} a_{\bfj} W_{\bfj},\qquad \mJ^d=\{\bfj\in\mJ:\; \uu_{\bfj}\subset \{1,\ldots,d\}\},
$$
which can be interpreted as function spaces over the $d$-dimensional tensor product domain $X^d$.
Since $H^d_{\bfW,\bfa}$ can be identified with $H^\infty_{\bfW,\bfa^d}$, where $\bfa^d_{\bfj}=\bfa_{\bfj}$ for 
all $\bfj$ in $\mJ^d$ and $\bfa^d_{\bfj}=0$ otherwise, we automatically get
$$
n_\epsilon(H^d_{\bfW,\bfa},H^d_{\bfW,\bfb})=n_\epsilon(H^d_{\bfW,\bfa^d},H^d_{\bfW,\bfb^d})=|\mJ_{\bfc^d,\epsilon}\cap \mJ^d|\le
n_\epsilon(H^\infty_{\bfW,\bfa},H^\infty_{\bfW,\bfb})=|\mJ_{\bfc,\epsilon}|. 
$$
In other words, control of the $\epsilon$-dimension for function spaces on $X^\infty$ guarantees automatically control for their restrictions to $X^d$,
simultaneously for all $d$. 

There is also a simple converse result which is in line with the main motivating observation that triggered this investigation, namely that dimension-robust results (e.g., estimates of $\epsilon$-dimensions as $\epsilon \to 0$ with constants that do not depend on $d$) for
$d$-dimensional approximation are possible if and only if a corresponding result holds for the infinite-dimensional approximation problem. To formulate it, we allow arbitrary rate measuring functions $\phi(\epsilon):\,(0,1] \to (0,\infty)$ with the property $\phi(\epsilon)\to\infty$ as $t\to 0$, where $\epsilon_0\in (0,\infty)$ is given.

\begin{pro}\label{prop3} Assume that $H^\infty_{\bfW,\bfa}$ is compactly embedded into $H^\infty_{\bfW,\bfb}$, where the weight sequence $\bfb$ is strictly positive (i.e., $b_{\bfj}>0$). Then, given any rate measuring function $\phi$, we have
$$
n^\infty_\epsilon(H^d_{\bfW,\bfa},H^d_{\bfW,\bfb})\le \phi(\epsilon),\qquad \epsilon\in (0,\epsilon_0],
$$
for all $d\in \mathbb{N}$, if and only if
$$
n^\infty_\epsilon(H^\infty_{\bfW,\bfa},H^\infty_{\bfW,\bfb})\le \phi(\epsilon),\qquad \epsilon\in (0,\epsilon_0].
$$
Similar statements hold if the requirement $\le \phi(\epsilon)$ for $\epsilon\in (0,\epsilon_0]$ is replaced by $=\mathrm{O}(\phi(\epsilon))$ or $=\mathrm{o}(\phi(\epsilon))$ for $\epsilon\to 0$.
\end{pro}

{\bf Proof}. The {\it if} direction is obvious, see above. For the {\it only if} direction, observe that, for fixed $\epsilon \in (0,\epsilon_0]$, the set $\mJ_{\bfc,\epsilon}$ must be finite, and thus
$\mJ_{\bfc,\epsilon}\subset \mJ^d$ for all $d\ge d_0$ with some finite $d_0=d_0(\epsilon)\in\mathbb{N}$ depending on $\epsilon$. Thus,
$$
\mJ_{\bfc,\epsilon}=\mJ_{\bfc,\epsilon}\cap \mJ^d=\mJ_{\bfc^d,\epsilon},\qquad d\ge d_0,
$$
and consequently
$$
n^\infty_\epsilon(H^\infty_{\bfW,\bfa},H^\infty_{\bfW,\bfb})\le n^\infty_\epsilon(H^{d_0}_{\bfW,\bfa},H^{d_0}_{\bfW,\bfb})\le \phi(\epsilon).
$$
Since $\epsilon\in (0,\epsilon_0]$ was fixed arbitrarily, this finishes the argument. 

\medskip
Clearly, Proposition \ref{prop3} is a statement about the asymptotics as $d\to \infty$. For practical reasons there is also valid interest in getting improved estimates for $n^\infty_\epsilon(H^{d}_{\bfW,\bfa},H^{d}_{\bfW,\bfb})$  in the pre-asymptotic range $d<d_0(\epsilon)$, and for situations, where $H^\infty_{\bfW,\bfa}$ is not compactly embedded into $H^\infty_{\bfW,\bfb}$, see, e.g., \cite{KSU}.

\medskip
We now turn to the study of a particular family of spaces based on redundant subspace splittings motivated by the study of multiscale approximation schemes. Assume that each $H_k$ possesses its own increasing ladder of finite-dimensional subspaces 
$$
V_{0,k}=\mathrm{span}\{1\} \subset V_{1,k}\subset \ldots\subset V_{j,k}\subset \ldots \subset H_k,\qquad j\in\mJ_k,
$$
and define the $H_k$-orthogonal complement spaces $W_{0,k}=V_{0,k}$, $W_{j,k}=V_{j,k}\ominus V_{j-1,k}$ (we also allow for finite $\mJ_k$ in which case we silently assume that $V_{K_k}=H_k$). In other words, we have $H_k=\bigoplus_{j\in\mJ_k} W_{j,k}$ with norm
$$
\|f\|_{H_k}^2 = \sum_{j=0}^\infty \|\ww_{j,k}\|_{H_k}^2,\qquad f\in H_k,
$$
where $f=\sum_{j=0}^{\infty} \ww_{j,k}$ is the unique $H_k$-orthogonal decomposition of $f$ with respect to the orthogonal subspaces $W_{j,k}$.
With the two systems $\{W_{j,k}\}$ and $\{V_{j,k}\}$ at hand, we can introduce two scales of Hilbert spaces, $H^\infty_{\bfW,\bfa}$ and 
$H^\infty_{\bfV,\bfa}$ as before. Recall from Section \ref{sec1} that, because of the redundancy of series representations with respect to $\bfV$, 
we started from the semi-norm
\be\label{Vsemi}
|f|^2_{\bfV,\bfa} := \inf_{\vv_{\bfj}\in V_{\bfj}:\; f = \sum_{\bfj\in \mJ} \vv_{\bfj}} \sum_{\bfj\in\mJ} a_{\bfj}\|\vv_{\bfj}\|_{\omega_{\bfj}}^2,
\ee
defined for $f\in\mathrm{span}\bfV$, with the infimum in (\ref{Vsemi}) taken with respect to all possible finite sum representations of $f$.
Conditions under which $|\cdot|_{\bfV,\bfa}$ is a norm are given below, together with an answer to the following question:
When are two spaces $H^\infty_{\bfV,\bfa}$ and $H^\infty_{\bfW,\bfb}$ equivalent, and what are the constants in the norm equivalences? 

We call a subset $\mJ'\subset\mJ$ monotone if $\bfi\in\mJ'$ implies $\bfj\in\mJ'$ for all $\bfj\le \bfi$, and we call a weight sequence $\bfa$ monotonously supported if 
its support $\omega_{\bfa}\subset \mJ$ is monotone. In the context of approximation processes with respect to $\bfV$, this is a natural assumption on the weight sequence.
Since $W_{\bfj}\subset V_{\bfj}$ for all $\bfj\in\mJ$, with the same norm, we always have
$$
|f|_{\bfV,\bfa}\le \|f\|_{\bfW,\bfa},
$$
which implies $H^\infty_{\bfW,\bfa}\subset H^\infty_{\bfV,\bfa}$ if the latter space is well-defined.
The opposite direction is more subtle (and often not true).
 
\begin{pro}\label{prop4} Let $\bfa$ and $\bfb$ be monotonously supported weight sequences with the same non-empty support.\\
a) The space $H^\infty_{\bfV,\bfa}$ is well-defined, i.e., (\ref{Vsemi}) defines a norm on $\mathrm{span}\bfV$, if and only if 
\be\label{NecCond}
\sum_{\bfi\in \omega_{\bfa}} a_{\bfi}^{-1}<\infty.
\ee
In particular, if (\ref{NecCond}) is violated then $|1|_{\bfV,\bfa}=0$. \\
b) Under the condition (\ref{NecCond}), we have the identity
\be\label{VequalW}
H^\infty_{\bfV,\bfa} = H^\infty_{\bfW,\hat{\bfa}}, \qquad \omega_{\hat{\bfa}}=\omega_{\bfa},\quad 
\hat{a}_{\bfj}:=(\sum_{\bfi\in\omega_{\bfa}:\,\bfi\ge \bfj} a_{\bfi}^{-1})^{-1},
\ee
with identical scalar products. 
\end{pro}

{\bf Proof}. Assume that (\ref{NecCond}) holds. For given $f\in \mathrm{span}\bfV$, take any of its finite representations
$$
f=\sum_{\bfj\in \mJ'} \vv_{\bfj},\qquad \mJ'\subset \omega_{\bfa},\quad |\mJ'|<\infty,
$$
with respect to $\bfV$. Each of the $\vv_{\bfj}\in V_{\bfj}$, $\bfj\in \mJ'$, possesses a finite decomposition
$$
\vv_{\bfj}=\sum_{\bfi\le\bfj} \ww_{\bfi,\bfj}:\quad \|\vv_{\bfj}\|^2_{\omega_{\bfj}}=\sum_{\bfi\le\bfj} \|\ww_{\bfi,\bfj}\|_{\omega_{\bfi}}^2,
$$
with respect to $\bfW$, where the $\ww_{\bfi,\bfj}\in W_{\bfi}$ are uniquely defined by construction. 
Then
$$
\sum_{\bfj\in \mJ'} a_{\bfj}\|\vv_{\bfj}\|_{\omega_{\bfj}}^2 =\sum_{\bfj\in \mJ'} a_{\bfj}\sum_{\bfi\le \bfj}\|\ww_{\bfi,\bfj}\|_{\omega_{\bfi}}^2
 = \sum_{\bfi\in \mJ''} \sum_{\bfj\in \mJ':\,\bfj\ge \bfi} a_{\bfj}\|\ww_{\bfi,\bfj}\|_{\omega_{\bfi}}^2,
$$
where $\mJ''\supset \mJ'$ is the smallest monotone subset of $\omega_{\bfa}$ containing $\mJ'$. Denote
$$
\ww_{\bfi}:=\sum_{\bfj\in \mJ':\,\bfj\ge \bfi} \ww_{\bfi,\bfj},\qquad \bfi\in \mJ''.
$$
Obviously, these $\ww_{\bfi}$ provide the unique decomposition 
$$
f=\sum_{\bfi\in \mJ''} \ww_{\bfi}
$$
of $f$ with respect to $\bfW$. Thus,
\begin{align}
|f|_{\bfV,\bfa}^2&= \label{hilf0}\inf_{\mJ': \,f=\sum_{\bfj\in \mJ'} \vv_{\bfj}} \sum_{\bfi\in \mJ''} \sum_{\bfj\in \mJ':\,\bfj\ge \bfi} a_{\bfj}\|\ww_{\bfi,\bfj}\|_{\omega_{\bfi}}^2 \\
&= \sum_{\bfi\in \omega_{\bfa}} (\inf_{\ww'_{\bfi,\bfj}\in W_{\bfj},\,\bfj\ge \bfi:\; \ww_{\bfi} =\sum_{\bfj\ge\bfi} \ww'_{\bfi,\bfj}} \sum_{\bfj\ge\bfi} a_{\bfj}\|\ww'_{\bfi,\bfj}\|_{\omega_{\bfi}}^2) \nonumber\\
&= \label{hilf} \sum_{\bfi\in \omega_{\bfa}} (\sum_{\bfj\in \omega_{\bfa}:\,\bfj\ge \bfi} a_{\bfj}^{-1})^{-1} \|\ww_{\bfi}\|_{\omega_{\bfi}}^2 = 
\sum_{\bfi\in \omega_{\bfa}} \hat{a}_{\bfi}\|\ww_{\bfi}\|_{\omega_{\bfi}}^2 = \|f\|_{\bfW,\hat{\bfa}}^2,   
\end{align}
where all summations are finite and the notation from (\ref{VequalW}) is used. 

The switch from the infimum for arbitrary $\bfV$ decompositions of $f$ in (\ref{hilf0}) to the sum of infimums for arbitrary decompositions
of the unique $\ww_{\bfi}$ into elements $\ww'_{\bfi,\bfj}$ in $W_{\bfi}$ is justified, since any decomposition of the latter type produces also a
finite decomposition of $f$ of the former type, i.e.,
$$
f=\sum_{\bfj\in \omega_{\bfa}} \vv'_{\bfj},\qquad \vv'_{\bfj} =\sum_{\bfi\le\bfj} \ww'_{\bfi,\bfj}.
$$
The last step leading to (\ref{hilf}) follows from the fact that in any Hilbert space  $H$
\be\label{Fact}
\inf_{u_j\in H:\;u=\sum_j u_j} \sum_{j} a_j\|u_j\|_H^2 = \left\{ \ba{ll} 0,& \sum_j a_j^{-1}=\infty,\\
&\\
(\sum_j a_j^{-1})^{-1}\|u\|_{H}^2, & \sum_j a_j^{-1}<\infty,\ea\right.\qquad u\in H,
\ee
see \cite[Lemma 3.1]{GHO}, which we apply with $u=\ww_{\bfi}$ and $H=W_{\bfi}$ for the finitely many $\bfi$ with $\ww_{\bfi}\neq 0$ (recall that
$f\in\mathrm{span}\bfV=\mathrm{span}\bfW$). Since (\ref{NecCond}) was assumed to hold, all weights $\hat{a}_{\bfj}$ in (\ref{VequalW}) are well-defined,
and the second case in (\ref{Fact}) gives the equality in (\ref{hilf}). Since $\|\cdot\|_{\bfW,\hat{\bfa}}$ is a norm, we find that, under the condition
(\ref{NecCond}), $|\cdot|_{\bfV,{\bfa}}$ is also a norm. This proves the sufficiency in Part a). From the shown equality of norms, Part b)
follows as well.

If (\ref{NecCond}) does not hold, we take any finite monotone index set $\mJ'\subset \omega_{\bfa}$ and decompose
$$
1 = \sum_{\bfj\in\mJ'} \vv_{\bfj},\qquad \vv_{\bfj} = a_{\bfj}^{-1}(\sum_{\bfi\in \mJ'} a_{\bfi}^{-1})^{-1},\quad \bfj\in\mJ'.
$$
This is an admissible finite decomposition of $1\in V_{\mathbf{0}}$ into multiples of $1$ belonging to the respective $V_{\bfj}$. Since by definition 
$\|1\|_{\omega_{\bfj}}=1$ for all $\bfj$, from this decomposition we have
$$
|1|_{\bfV,\bfa}^2\le \sum_{\bfj\in \mJ'} a_{\bfj}\|\vv_{\bfj}\|_{\omega_{\bfj}}^2 = \sum_{\bfj\in \mJ'} a_{\bfj}^{-1}(\sum_{\bfi\in \mJ'} a_{\bfi}^{-1})^{-2}
=(\sum_{\bfi\in \mJ'} a_{\bfi}^{-1})^{-1}.
$$
Thus, letting $\mJ'$ grow, we see that $|1|_{\bfV,\bfa}=0$, which shows that $|\cdot|_{\bfV,\bfa}$ is not a norm. This completes the proof of
Proposition \ref{prop4}.

\medskip
As an immediate consequence of Propositions \ref{prop1} and \ref{prop4} we obtain necessary and sufficient conditions for embeddings of $H^\infty_{\bfV,\bfa}$ or $H^\infty_{\bfW,\bfa}$ into $H^\infty_{\bfV,\bfb}$ or $H^\infty_{\bfW,\bfb}$.  
\begin{cor}\label{cor1} We have $H^\infty_{\bfV,\bfa}\subset H^\infty_{\bfW,\bfb}$
with the embedding operator bounded by $C<\infty$ if and only if
$$
\sum_{\bfi\le \bfj\in \mJ} a_{\bfj}^{-1}\le C^2 b_{\bfi}^{-1},\qquad \forall\,\bfi\in \omega_{\bfa}\subset \omega_{\bfb}.
$$
In particular, $H^\infty_{\bfV,\bfa}= H^\infty_{\bfW,\bfa}$ with equivalent norms if and only if
\be\label{VaWa}
\sum_{\bfj\le \bfi\in \mJ} a_{\bfi}^{-1}\le C^2 a_{\bfj}^{-1},
\ee
for some $C$. Then $C^2$ is a bound for the relative spectral condition number of the scalar products given
by the two norms $\|\cdot\|_{\bfV,\bfa}$ and $\|\cdot\|_{\bfW,\bfa}$ on $H^\infty_{\bfV,\bfa}=H^\infty_{\bfW,\bfa}$.
\end{cor}

\medskip
The condition (\ref{VaWa}) is very stringent and holds only for special classes of weights $\bfa$, see Section \ref{sec4}.
The result from (\ref{VequalW}) also shows that defining equivalent splitting norms for natural subspaces of $H^\infty_{\bfV,\bfa}$ is a subtle issue.
For example, take any finite monotone subset $\mJ'$ of $\mJ$ and consider
$$
H^\infty_{\bfV,\bfa}|_{\mJ'} := \{ f\in H^\infty_{\bfV,\bfa}: \; \vv_{\bfj}=0\mbox{ if $\bfj\not\in \mJ'$}\}.
$$
Can this subspace be identified with some $H^\infty_{\bfV,\bfa'}$, where the new weight sequence $\bfa'$ is supported in $\mJ'$?
We would be happy with a characterization with equivalent norms, where the constants in the norm equivalence 
$$
\|f\|_{\bfV,\bfa}\approx \|f\|_{\bfV,\bfa'}, \qquad f\in H^\infty_{\bfV,\bfa}|_{\mJ'},
$$
do not depend on $\mJ'$ and $\bfa$.
Unfortunately, this is not always possible, see \cite{GHO} for a discussion in a special case. Indeed, since we have (\ref{VequalW}) and in view of Proposition \ref{prop2} (now applied to $H^\infty_{\bfV,\bfa'}$), any weight sequence
$\bfa'$ providing the desired norm equivalence must satisfy
$$
\sum_{\bfj\in \mJ':\, \bfj\ge \bfi} (a'_{\bfj})^{-1} \approx \hat{a}_{\bfi}^{-1}=\sum_{\bfj\in \mJ: \bfj\ge \bfi} a_{\bfj}^{-1},
$$
uniformly for all $\bfi\in\mJ'=\mathrm{supp}(\bfa')$, with constants independent of  $\mJ'$ and $\bfa$. The proof of \cite[Theorem 3.4]{GHO} provides 
a linear programming approach to check this condition, and a set of weights for which uniform and $\mJ'$ independent bounds in the norm equivalence
cannot be achieved.

\section{Case Studies}\label{sec3}

In this section, we give illustrative examples of how to put the developed machinery to work. We first consider decomposition spaces based on tensor product splines and apply the results of Section \ref{sec2} to derive embedding conditions and to estimate $\epsilon$-dimensions. Then we give an example of comparing $H_{\bfW,\bfa}$ and $H_{\bfW',\bfa}$ norms for slightly different 
subspace families $\bfW$ and $\bfW'$, sharing the same underlying sequence of coordinate Hilbert spaces $H_k$. The question was triggered by work in \cite{HR,HRW,HS} devoted to the equivalence of ANOVA and anchored decompositions of 
Sobolev spaces of mixed smoothness for infinitely many variables. Finally, the approximation of maps between infinite-dimensional spaces
is discussed for the example of least-squares regression for Hilbert space valued mappings. 

\subsection{Tensor product spline decompositions}\label{sec4}

Let us examine the case when the coordinate spaces $H_k$ coincide with $L_2$-spaces over the interval $X_k=I=[0,1]$ (periodic or non-periodic),
and the spaces $V_{j,k}$ are spanned by splines of polynomial degree $p_k\in \mathbb{Z}_+$ and smoothness $r_k\le p_k-1$ on uniform dyadic partitions of step-size $2^{-j}$, $j\ge 0$. Choosing $r_k=-1$ corresponds to the case of non-smooth piecewise polynomial functions, $r_k=p_k -1$ to smooth splines, while $r_k=0$ and $r_k=1$ covers finite element applications to second- and fourth-order elliptic boundary value problems. To conform with previous assumptions on $V_{0,k}=W_{0,k}$, we need to modify $V_{0,k}$ and $V_{1,k}$ such that $V_{0,k}$ just contains constant functions. The case when $V_{0,k}=\mathrm{span}\{1\}$ automatically holds is the case of smooth periodic splines ($r_k=p_k-1\ge 0$) which we will from now on concentrate on.
The extension to general periodic and non-periodic spline spaces as well as to spline spaces on quasi-uniform partitions  is more or less straightforward. 

The spaces $V_{j,k}$ can be equipped with locally supported stable B-spline bases, and the $L_2$ orthogonal complement spaces 
$W_{j,k}$ also possess locally supported, so-called prewavelet basis functions. It is easy to check that 
$$
\dim V_{j,k}= 2\dim W_{j,k}= 2^j,\quad j\ge 1, \qquad \dim V_{0,k}= \dim W_{0,k} =1.
$$
The dimension formulas are well-known also for other families of univariate spline spaces but are more involved. Up to additive constants, they are of the order $(p_k-r_k)2^j$.
The following lemma characterizes univariate Sobolev spaces, and is instrumental for many applications. In its formulation, we temporarily drop the 
subscript $k$.
\begin{lem}\label{lem2} Let $\{W_j\}_{j\ge 0}$ be the $L_2$ orthogonal subspace family associated with smooth periodic splines of degree $p\ge 1$ over uniform dyadic partitions of the unit interval $I$, as introduced above.\\
For $|s|<p+1/2 $, there are positive constants $0<\underline{\lambda}_{s,p}\le \bar{\lambda}_{s,p}<\infty$ such that, for any $u\in H^s(I)$, there is a unique  representation $u= \sum_{j=0}^\infty \ww_j$ (convergence in the sense of $L_2(I)$ if $s\ge 0$, otherwise
in a distributional sense) with $\ww_j\in W_j$ for $j\ge 0$ such that
$$
\|\ww_0\|_{L_2}^2 + \underline{\lambda}_{s,p}\sum_{j=1}^\infty 2^{2sj}\|\ww_j\|_{L_2}^2\le \|u\|^2_{H^s(I)} \le \|\ww_0\|_{L_2}^2 + \bar{\lambda}_{s,p}\sum_{j=1}^\infty 2^{2sj}\|\ww_j\|_{L_2}^2.
$$
If $s=0$, we have by construction the equality 
$$
\|u\|^2_{L_2}= \|\ww_0\|_{L_2}^2 + \sum_{j=1}^\infty \|\ww_j\|_{L_2}^2  \qquad (\underline{\lambda}_{0,p}=\bar{\lambda}_{0,p}=1).
$$
\end{lem}
This lemma allows us to relate Sobolev spaces of mixed dominating smoothness on the cubes $X^d$ and $X^\infty$   to
abstract $H^d_{\bfW,\bfa}$ and $H^\infty_{\bfW,\bfa}$ spaces if $\bfW$ is generated from subspace families $\{W_{j,k}\}_{j\ge 0}$ of 
smooth periodic splines of
degree $p_k$, $k\in \mathbb{N}$. 
The norm equivalence of Lemma \ref{lem2} suggests 
to look at weight sequences $\bfa$ given by $a_{\mathbf{0}}=1$ and
\be\label{Aweights}
a_{\bfj} = \gamma_{\omega_{\bfj}}^{-1} \prod_{k\in \omega_{\bfj}} \underline{\lambda}_{s_k,p_k}2^{2s_kj_k},\qquad \mathbf{0}\neq \bfj\in\mJ,
\ee
where the sequence $\gamma:=\{\gamma_\omega\}_{\omega\in\mP_f(\mathbf{N})}$ is non-negative. If we are in the popular product-weight case
\be\label{prodweights}
\gamma_{\omega}=\prod_{k\in \omega} \gamma_k,\quad \emptyset\neq \omega\subset \mP_f(\mathbb{N}), 
\ee
then the weighted mixed-norm Sobolev space 
$$
H^{\bfs}_\gamma(X^\infty):=\bigotimes_{k=1}^\infty \gamma_k H^{s_k}(X_k),\qquad \bfs=(s_1,s_2,\ldots)\ge 0,\; \gamma=(\gamma_1,\gamma_2,\ldots)>0,\; X_k=I,
$$
is norm-$1$ embedded into $H_{\bfW,\bfa}(X^\infty)$. Here, slightly abusing notation, we interpret $\gamma H^{s}(I)$ as $H^{s}(I)$ equipped with the
norm 
$$
\|u\|_{\gamma H^{s}(I)}^2 = \|P_0u\|_{L_2}^2 + \gamma^{-1}|u-P_0u|_{H^s}^2,
$$
where $P_0$ is the $L_2$ orthoprojector onto $\mathrm{span}\{1\}$. Thus, $\epsilon$-dimensions for approximation in $L_2$ or other $H^{\tilde{\bfs}}_{\tilde{\gamma}}(X^\infty)$ norms of elements $H^{\bfs}_\gamma(X^\infty)$ 
can be estimated via Proposition \ref{prop2}. For $\gamma_k=1$, and under special conditions on $s_k\to \infty$ (penalization by increasing smoothness
as $k\to \infty$), this question has been examined in 
\cite{DG} in the similar case of Fourier decompositions. 

Here we examine the complementary situation when the smoothness $s_k=s>0$ is fixed, i.e., $\bfs=(s,s,\ldots)$. and the penalization is realized by only assuming a sufficient decay of the weights $\gamma_\omega$ in (\ref{Aweights}). To simplify notation, we also fix the spline degree
$p_k=p$, and set $\underline{\lambda}:=\underline{\lambda}_{s,p}$).   
Note that in this case $\gamma_\omega\to 0$ (meaning that for any $\delta>0$ the inequality $\gamma_\omega\ge \delta$
holds for only a finite number of $\omega\in\mP(\mathbb{N})$) is necessary and sufficient for the compact embedding of $H^{\bfs}_\gamma(X^\infty)$ into $L_2(X^\infty)$. From the definition of $\mJ_\epsilon$ and since
$$
\dim W_{\bfj} = \prod_{k\in\omega_{\bfj}} \dim W_{j_k} =2^{|\bfj|_1 -|\bfj|_0},\qquad |\bfj|_1:=\sum_{k\in \omega_{\bfj}} j_k,\qquad \mathbf{0}\neq\bfj\in\mJ,
$$
we find that
\bea
n_\epsilon(H^\infty_{\bfW,\bfa},L_2(X^\infty))&=&1+\sum_{\mathbf{0}\neq \bfj \in \mJ_\epsilon} 2^{|\bfj|_1-|\bfj|_0}= 
1+\sum_{\emptyset\neq \omega\in \mP_f(\mathbb{N})} \left(\sum_{\bfj\in\mJ_\epsilon:\;\omega_{\bfj}=\omega} 2^{|\bfj|_1-|\bfj|_0}\right)\\
&=&1+\sum_{\emptyset\neq \omega\in \mP_f(\mathbb{N})}\; \left(\sum_{\bfj\in\mJ:\;\omega_{\bfj}=\omega\; 2^{2s|\bfj|_1 } 
\le \epsilon^{-2}\gamma_\omega \underline{\lambda}^{-|\omega|}} 2^{|\bfj|_1-|\bfj|_0}\right)\\
&=&1+\sum_{\emptyset\neq \omega\in \mP_f(\mathbb{N})}\; \left(\sum_{\bfj\in\mJ:\;\omega_{\bfj}=\omega,\; 2^{|\bfj|_1 -|\bfj|_0} 
\le \epsilon^{-1/s}(\gamma_\omega (2\underline{\lambda})^{-|\omega|})^{1/(2s)}} 2^{|\bfj|_1-|\bfj|_0}\right)\\
&=&1+\sum_{\emptyset\neq \omega\in \mP_f(\mathbb{N}): m(\omega)\ge 0}\;\left(\sum_{m=0}^{m(\omega)}  2^{m}|\omega|^m\right)\\
&\le& 1+ 2\sum_{\emptyset\neq \omega\in \mP_f(\mathbb{N}): m(\omega)\ge 0} (2|\omega|)^{m(\omega)},
\eea
where the integer $m(\omega)$ is defined by
$$
m(\omega):=[\log_2(\epsilon^{-1/s}(\gamma_\omega(2\underline{\lambda})^{-|\omega|})^{1/(2s)})].
$$
We have used that there are exactly $|\omega|^m$ index vectors $\bfj$ with support $\omega=\omega_{\bfj}$ and the same value $m=|\bfj|_1-|\bfj|_0$. Note that we get a lower bound for $n_\epsilon(H^\infty_{\bfW,\bfa},L_2(X^\infty))$ if the forefactor $2$ in front of the last sum is replaced by $1$.
In other words, we still have an optimal bound (within a small absolute factor). 

Obviously $n_\epsilon(H^\infty_{\bfW,\bfa},L_2(X^\infty))$ is finite if and only if $m(\omega)<0$ for all but finitely many $\omega\in \mP_f(\mathbb{N})$, i.e., if $\gamma_\omega\to 0$. To give quantitative estimates for $\epsilon$-dimensions, we need additional assumptions on the weights $\gamma_\omega$. 
Let us start with the case where we assume $\gamma_\omega=0$ for all $\omega$ with $|\omega|>1$. 
Then, for $\epsilon \to 0$, by substituting the formula for $m(\omega)$ we have 
$$
n_\epsilon(H^\infty_{\bfW,\bfa},L_2(X^\infty))\le 1+ \sum_{\omega: |\omega|=1,\,m(\omega)\ge 0} 2^{m(\omega)}
\approx 1 + \epsilon^{-1/s} (2\underline{\lambda})^{-1/(2s)} \sum_{k:\, \gamma_{\{k\}}>2\underline{\lambda}\epsilon^2}
\gamma_{\{k\}}^{1/(2s)}.
$$
Thus, the best possible rate $\mathrm{O}(\epsilon^{-1/s})$ is achievable if and only if the sequence $\gamma_{\{k\}}^{1/(2s)}$, $k\in\mathbb{N}$, is summable.

Weaker assumptions lead to a deterioration of the growth rate of the $\epsilon$-dimension.
Already in the case when $\gamma_\omega=0$ holds only for $|\omega|>2$, we obtain
$$
n_\epsilon(H^\infty_{\bfW,\bfa},L_2(X^\infty))\approx 1 + (2\underline{\lambda}\epsilon^2)^{-1/(2s)} \sum_{k:\, \gamma_{\{k\}}>2\underline{\lambda}\epsilon^2}
\tga_{\{k\}}^{1/(2s)} + (2\underline{\lambda}\epsilon^2)^{-1/s}\sum_{k\neq l:\, \gamma_{\{k,l\}}>2\underline{\lambda}\epsilon^2}
\tga_{\{k,l\}}^{1/s}.
$$
I.e., if there is at least one positive weight $\gamma_{\{k,l\}}>0$ then the growth rate of $\epsilon$-dimensions asymptotically increases to the order of $\epsilon^{-2/s}$.
Similarly, for each $\gamma_\omega>0$ we always have a lower bound of 
$$
n_\epsilon(H^\infty_{\bfW,\bfa},L_2(X^\infty))\ge \tga_{\omega}^{(1+\log_2 |\omega|)/(2s)} \epsilon^{-(1+\log_2 |\omega|)/s}, \qquad c_\omega>0,
$$
for small enough $0<\epsilon\le \epsilon_\omega$. Thus, no rate estimate of the form $\mathrm{O}(\epsilon^{-\alpha})$ with finite $\alpha$ can be expected if
$\gamma_\omega>0$ for $\omega\in \mP_f(\mathbb{N})$ of arbitrarily large size $|\omega|$. We leave it to the reader to discuss non-polynomial growth estimates, e.g., for the case of
summable product weights, where one can obtain bounds of the form
$$
n_\epsilon(H^\infty_{\bfW,\bfa},L_2(X^\infty)) =\mathrm{O}( (C\log_2 (1/\epsilon))^{\log_2 (1/\epsilon)} ), \quad \epsilon\to 0,
$$
involving a weight-dependent constant $C$. To summarize, compact embedding can be achieved by introducing weights damping the influence of certain coordinate directions (as in the case of product 
weights), or by penalizing certain coordinate combinations (weights that depend on $\omega\in \mP_f(\mathbb{N})$ in a more general way, see \cite{Gu,Wa1}
for examples). An alternative is to penalize certain coordinate directions by increasing the smoothness parameters $s_k\to \infty$, as done in \cite{DG,GH}. 

\medskip
We now turn to decomposition spaces related to $\bfV$ that are generated by the subspace families $\{V_{j,k}\}_{j\ge 0}$, $k\in\mathbb{N}$ which, despite
their redundancy, often lead to simpler approximation algorithms for $d$-dimensional problems. We examine the conditions
under which $H^\infty_{\bfV,\bfa}$ is well-defined for the weights given by (\ref{Aweights}), and when $H^\infty_{\bfV,\bfa}$ is  the same space (up to equivalent norms) as 
$H^\infty_{\bfW,\bfa}$. To simplify notation,
set $\tga_\omega:=\gamma_\omega \prod_{k\in\omega} \underline{\lambda}_{s_k,p_k}^{-1}$ for $\omega\in \mP_f(\mathbb{N})$.
Checking (\ref{NecCond}) in Proposition \ref{prop4} reveals that for positive weight sequences $\gamma>0$ the condition
$$
\sum_{\bfj\in\mJ} a_{\bfj}^{-1}=\sum_{\bfj\in\mJ} \tga_{\omega_{\bfj}}2^{-2\bfs\cdot\bfj} = \sum_{\omega\in\mP_f(\mathbb{N})} \tga_\omega \sum_{\bfj:\,\omega_{\bfj}=\omega} 2^{-2\bfs\cdot\bfj}
=\sum_{\omega\in\mP_f(\mathbb{N})} \tga_\omega \prod_{k\in \omega}\frac{2^{-2s_k}}{1-2^{-2s_k}} < \infty
$$
guarantees that $H^\infty_{\bfV,\bfa}$ is well-defined. Moreover, under this condition, we have $H_{\bfV,\bfa}^\infty = H_{\bfW,\hat{\bfa}}^\infty$ with the weight sequence $\hat{\bfa}$ defined by
$$
\hat{a}_{\bfi}^{-1} := \sum_{\bfj\le \bfi} \tga_{\omega_{\bfj}}2^{-2\bfs\cdot\bfj} =
\sum_{\omega\in \mP_f(\mathbb{N}):\,\omega_{\bfi}\subset\omega } \tga_\omega \sum_{\bfj\ge \bfi:\,\omega_{\bfj}=\omega} 2^{-2\bfs\cdot\bfj}
=2^{-2\bfs\cdot\bfi} \sum_{\omega\in \mP_f(\mathbb{N}): \omega_{\bfi}\subset \omega} \tga_\omega \prod_{k\in \omega}\frac{1}{1-2^{-2s_k}}.
$$
According to (\ref{VaWa}) in Proposition \ref{prop4}, since
$$
a_{\bfi}\hat{a}_{\bfi}^{-1}=\tga_{\omega_{\bfi}}^{-1}2^{2\bfs\cdot\bfi}\hat{a}_{\bfi}^{-1}=\sum_{\omega\in \mP_f(\mathbb{N}):\,\omega_{\bfi}\subset\omega} 
\frac{\tga_\omega}{\tga_{\omega_{\bfi}}} \prod_{k\in \omega}\frac{1}{1-2^{-2s_k}}\ge
\prod_{k\in \omega_{\bfi}}\frac{1}{1-2^{-2s_k}},
$$
for $H^\infty_{\bfV,\bfa}=H^{\infty}_{\bfW,\bfa}$ to hold (with equivalent norms) it is necessary that
\be\label{Sgrowth}
\sum_{k=1}^\infty 2^{-2s_k} <\infty.
\ee
If, in addition,
$$
\sum_{\omega\in \mP_f(\mathbb{N}): \omega_{\bfi}\subset \omega} \tga_\omega \le C \tga_{\omega_{\bfi}},
$$
then we arrive at a sufficient condition. In particular, for product-weights (\ref{prodweights}) satisfying the summability condition
$$
\sum_{k=1}^\infty \gamma_k \underline{\lambda}_{s_k,p_k}^{-1} < \infty,
$$
the growth condition (\ref{Sgrowth}) on the assumed coordinate-wise smoothness exponents $s_k$ is necessary and sufficient for 
$H^\infty_{\bfV,\bfa}=H^{\infty}_{\bfW,\bfa}$.

Estimates of $\epsilon$-dimensions and norm equivalences can be obtained in a similar way in other situations. For instance, motivated by the definition of anisotropic Sobolev spaces, instead of mixed-type norms, one could consider weighted norms of the form
\be\label{DefWIsotrop}
\|f\|^2_{\bfW,\bfa} = \sum_{\bfj\in\mJ} a_{\bfj}\|\ww_{\bf j}\|_{L_2}^2, \qquad  a_{\bfj}= \gamma_{\omega_{\bfj}}^{-1} \sum_{k\in \omega_{\bfj}} 2^{2s_k j_k},
\ee
where for $\bfW$ we can take one of the above subspace families. 
Another case of interest are decomposition spaces related to interpolation processes. To this end, define projectors $I_{j,k}: \, C(I)\to V_{j,k}$ by interpolation at the knots (for odd $p_k$) and at the interval midpoints (for even $p_k$) of the uniform partitions underlying the spline space
$V_{j,k}$, respectively, and set $\tW_{j,k}=\mathrm{Ran}(I_{j,k}-I_{j-1,k})$, $j\ge 1$, $\tW_{0,k}=\mathrm{span}\{1\}$. Then we have a counterpart of 
Lemma \ref{lem2} (we drop again the subscript $k$): For the range $1/2 < s < p+1/2$, the norm equivalence
$$
\|\tw_0\|_{L_2}^2 + \underline{\lambda}'_{s,p}\sum_{j=1}^\infty 2^{2sj}\|\tw_j\|_{L_2}^2\le \|u\|^2_{H^s(I)} \le \|\tw_0\|_{L_2}^2 + \bar{\lambda'}_{s,p}\sum_{j=1}^\infty 2^{2sj}\|\tw_j\|_{L_2}^2
$$
holds for all $u\in H^s(I)$ which are uniquely decomposed as
$$
u=\sum_{j=0}^\infty \tw_j,\qquad \tw_0=I_0u,\quad \tw_j=I_ju-I_{j-1}u\in \tW_{j},\quad j\ge 1.
$$
Since $s>1/2$ we have $H^s(I)\subset C(I)$, which guarantees the well-posedness of this decomposition. The constants $\underline{\lambda}'_{s,p},\bar{\lambda'}_{s,p}$ generally differ from those in Lemma \ref{lem2}. Thus, we can also define spaces $H^{\infty}_{\tilde{\bfW},\bfa}$, where $\tilde{\bfW}$ is induced by tensorization of the subspace families $\{W_{j,k}\}_{j\ge 0}$, $k\in\mathbb{N}$, and study approximation processes based on interpolation and function evaluation for functions of infinitely many variables.

\subsection{Equivalence of ANOVA and anchored spaces}\label{sec5} 

In this section, we consider a particular example from the literature where $H^\infty_{\bfW,\bfa}$ spaces have appeared in the analysis of
infinite-dimensional integration problems and quasi Monte Carlo methods. Although there is a general formulation in terms of weighted tensor product RKHS, see \cite{GMR,HR,KSWW}, we consider here only the case of first order mixed dominating smoothness when the $H_k$ coincide with the Sobolev space $H^1(I)$ on the same unit interval $X_k=I=[0,1]$, equipped with a particular scalar product.
We fix a projector $P:\,H^1(I)\to W_{0}:=\mathrm{span}\{1\}$ onto the subspace of constant functions, set $W_{1,k}=\mathrm{Ran}(I-P_k)$, and introduce a $P_k$-dependent norm on $H^1(I)$ by
$$
\|f\|_{P_k}^2 = |P_kf|^2 + |f-P_kf|_{H^1}^2,
$$
where $|f|^2_{H^1}=\int_I |f'(x)|^2\,dx$ is the standard semi-norm on $H^1(I)$.
Thus, we are in the setting of Section \ref{Def}, if we identify $H_k$ with $H^1(I)=W_{j,0}+W_{1,k}$ equipped with the norm $\|\cdot\|_{P_k}$, and set $\mJ_k=\{0,1\}$, $k\in \mathbb{N}$. Since $\mJ=\{0,1\}^\infty$ can be identified with $\mP_f(\mathbb{N})$, we will simplify notation and subscript $\bfW$ and $\bfa$
with finite subsets of $\mathbb{N}$, i.e.,
\be\label{WAgamma}
\bfW=\{W_\omega:=\bigotimes_{k\in\omega} W_{1,k}\}_{\omega\in \mP_f(\mathbb{N})},\qquad \bfa=\{a_\omega:=\gamma_\omega^{-1}\}_{\omega\in \mP_f(\mathbb{N})}.
\ee

The resulting spaces $H^\infty_{\bfW,\bfa}$ depend on the choice of the projectors $P_k$, 
which poses the question of comparing them for different $\{P_k\}_{k\in\mathbb{N}}$ and different weight sequences $\gamma$.  This question was addressed in \cite{HR} (and later in \cite{HRW,HS} for the $L_p$ setting) for the two most prominent projectors 
$$
P_kf:=\int_I f(t_k)\,dt_k,\qquad  P'_kf:=f(x^\ast),
$$
where $x^\ast\in I$ is fixed. 

For simplicity (and as in \cite{HR,HRW,HS}), we take the same projector for all directions. Denote by $P_{\omega}=\otimes_{k\in \omega} P_k$ and 
$(\mathrm{Id}-P)_{\omega}=\otimes_{k\in \omega} (\mathrm{Id}-P_k)$ the projectors for functions of the tensor product spaces
associated with the variable group indexed by $\omega\subset \mathbb{N}$, 
similarly for $P'_\omega$ and $(\mathrm{Id}-P')_\omega$. For coordinate vectors $x\in X^\infty$, we use the notation
$x=(x_{\omega_1},x_{\omega_2},\ldots,x_{\omega_m})$ to indicate that the coordinate directions are organized in groups of variables 
$x_{\omega_l}\in X^{\omega_l}$  with index sets $\omega_l$, $l=1,\ldots,m$, forming a partition
of $\mathbb{N}$. E.g., 
$$
P_{\omega^c}f(x_\omega)=\int_{X^{\omega^c}} f(x_\omega,t_{\omega^c})\,dt_{\omega^c},\qquad P'_{\omega^c}f(x_\omega)= f(x_\omega,x^\ast_{\omega^c}),
$$
for any finite index set $\omega\subset \mathbb{N}$. Here, $\omega^c:=\mathbb{N}\backslash \omega$ denotes the complementary set. 
These projectors define (under certain technical conditions) the ANOVA decomposition
$$
f(x) = \sum_{\omega\in \mP_f(\mathbb{N})} f_\omega(x_\omega),\qquad f_\omega:=(\mathrm{Id}-P)_\omega P_{\omega^c}f=\sum_{\omega'\subset \omega} (-1)^{|\omega|-|\omega'|} P_{\omega'}P_{\omega^c}f,
$$
and the anchored decomposition 
$$
f(x) = \sum_{\omega\in \mP_f(\mathbb{N})} f'_\omega(x_\omega),\qquad f'_\omega:=(\mathrm{Id}-P')_\omega P'_{\omega^c}f=\sum_{\omega'\subset \omega} 
(-1)^{|\omega|-|\omega'|} P'_{\omega'} P'_{\omega^c}f
$$
of functions of infinitely many variables. While the ANOVA decomposition has a straightforward justification in statistical terms, its evaluation  
requires integration which is more costly than evaluating function values as needed for anchored decompositions. This is the main motivation
for comparing properties of these decompositions.

Now, given any weight sequence $\gamma=\{\gamma_\omega\}$ with $\gamma_{\emptyset}=1$ and monotone support $\omega_\gamma:=\{\omega\in\mP_f(\mathbb{N}):\, \gamma_\omega>0\}$, we introduce weighted ANOVA and anchored decomposition norms by setting
\be\label{AanNorms}
\|f\|_{\gamma,A}^2:= \sum_{\omega} \gamma_\omega^{-1}\|f^{(\omega)}_\omega(x_\omega)\|_{L_2}^2,
\qquad \|f\|_{\gamma,an}^2:= \sum_{\omega} \gamma_\omega^{-1}\|{f'}^{(\omega)}_\omega(x_\omega)\|_{L_2}^2,
\ee
respectively. Here, $f_\omega^{(\omega)}=(P_{\omega^c}f)^{(\omega)}$ denotes the mixed first derivative with respect to all variables with indices from $\omega$, similarly for ${f'}_\omega^{(\omega)}$. These norms are special instances of $H^\infty_{\bfW,\bfa}$ norms, if the subspace family $\bfW$ and the weights $\bfa$ given by $\gamma$ as in (\ref{WAgamma}) are suitably defined. There are subtleties such as the dependence of $\mathcal{X}$ on $\gamma$, the pointwise meaning of the infinite summations, the practical meaning of $P_{\omega^c}$, etc., which we ignore here and instead refer to \cite{GMR,HR}. 

The question for which weight sequences $\gamma$ these norms define the same space has been investigated in \cite{HR}, and later in \cite{HRW,HS} for the $L_p$ setting, $1\le p\le\infty$. For $p=2$, necessary and sufficient conditions are established in \cite{HR} for special classes of weights
while in \cite{HS} a sufficient condition for general $\gamma$ is obtained by interpolating between the $L_1$ and $L_\infty$ results from \cite{HRW}. Our goal is to give a direct proof in the $L_2$ case.
\begin{pro}\label{prop5} For a given weight sequence  $\gamma=\{\gamma_\omega\}$ with monotone support $\omega_\gamma$, ANOVA and anchored norms are equivalent if there exist a sequence $\alpha=\{\alpha_\omega\}$ and constants $0<C',C''<\infty$
such that
\be\label{L2cond0}
\sum_{\homega\in\omega_\gamma:\,\homega\supset \omega} q^{|\homega|}\alpha_{\homega} \le C'q^{|\omega|}\alpha_{\omega},
\ee
and 
\be\label{L2cond}
\sum_{\omega'\in\omega_\gamma:\,\omega'\subset \omega} \frac{\alpha_{\omega'}}{\gamma_{\omega'}} \le C''\frac{\alpha_{\omega}}{\gamma_{\omega}},
\ee
for all $\omega\in\omega_\gamma$, where the constant $q=1/3-x^\ast (1-x^\ast)$ depends on the anchor $x^\ast$. If we set 
$\alpha_\omega=\tilde{q}^{|\omega|}\sqrt{\gamma_{\omega}}$ with some fixed $1\le  \tilde{q}\le 3/2$,
and $q=\frac13$ then (\ref{L2cond0}), (\ref{L2cond}) imply the sufficient conditions provided in \cite{HS} for $x^\ast=0$.
\end{pro}

{\bf Proof}. To avoid any discussions on pointwise definitions of infinite sums, we can silently assume that all considered $f$ have finite decompositions, and that all summations are only with respect to index sets contained in $\omega_\gamma$. We rely on the following representation formulas from Lemma 1 in \cite{HRW}: 
$$
f(x_\omega,x_{\omega^c})= \sum_{\omega'\subset \omega^c} \int_{X^{\omega'}} (P'_{(\omega'\cup \omega)^c}f)^{(\omega')}(x_\omega,t_{\omega'})\cdot \kappa_{\omega',an}(x_{\omega'},t_{\omega'})\,dt_{\omega'},\qquad \omega\subset \mP_f(\mathbb{N}),
$$
where $\kappa_{\omega',an}(x_{\omega'},t_{\omega'}) =\prod_{k\in \omega'} \kappa_{an}(t_k,x_k)$ are kernels defined by the identity
$$
g(x)=g(x^\ast) +\int_{I} \frac{dg}{dx}(t)\kappa_{an}(x,t)\,dt, \qquad \kappa_{an}(x,t)=\mathbf{1}_{[0,x]}(t)-\mathbf{1}_{[0,x^\ast]}(t), 
$$
where $g\in H^1(I)$ is arbitrary. Note that $\kappa_{an}$ is not a reproducing kernel, however, it is related to the reproducing kernel of $H^1(I)$ with norm induced by the anchored projector
via differentiation.

Similarly,
$$
f(x_\omega,x_{\omega^c})= \sum_{\omega'\subset \omega^c} \int_{X^{\omega'}} (P_{(\omega'\cup \omega)^c}f)^{(\omega')}(x_\omega,t_{\omega'})
\cdot \kappa_{\omega',A}(x_{\omega'},t_{\omega'})\,dt_{\omega'},
$$
where again $\kappa_{{\omega'},A}(x_{\omega'},t_{\omega'}) =\prod_{k\in {\omega'}} \kappa_{A}(t_k,x_k)$ are kernels, this time defined by
$$
g(x)=\int_I g(t)\,dt +\int_{I} \frac{dg}{dx}(t)\kappa_{A}(x,t)\,dt,\qquad \kappa_{A}(x,t)=\left\{\ba{ll} t,& 0\le t <x,\\ -(1-t),& x\le t\le 1.\ea \right.
$$

Let us show how to derive a sharp estimate of the ANOVA norm by the anchored norm, the other direction is completely analogous. Since $f^{({\omega})}_\omega(x_\omega)=(P_{\omega^c}f)^{(\omega)}(x_\omega)$,
we have from the above formula
$$
f^{(\omega)}_\omega(x_\omega)=\sum_{{\omega'}\subset \omega^c} \int_{X^{\omega'}} 
f^{(\omega\cup {\omega'})}(x_\omega,t_{\omega'},x^\ast_{(\omega\cup {\omega'})^c})\left(\int_{X^{\omega^c}} 
\kappa_{{\omega'},an}(x_{\omega'},t_{\omega'})\,dx_{\omega^c}\right)\,dt_{\omega'}.
$$
Then we check that 
$$
\int_{X^{\omega^c}} \kappa_{{\omega'},an}(x_{\omega'},t_{\omega'})\,dx_{\omega^c}
=\int_{X^{{\omega'}}} \kappa_{{\omega'},an}(x_{\omega'},t_{\omega'})\,dx_{\omega'} 
= \prod_{k\in {\omega'}}  K_{an}(t_k), \qquad \omega'\subset \omega^c,
$$
where 
$$
K_{an}(t):=\int_I \kappa_{an}(x,t)\,dx=\left\{\begin{array}{ll} -t, & 0\le t<x^\ast,\\ 1-t,& x^\ast\le t\le 1.\end{array} \right.
$$ 
Next, we define the constant $q$ entering the conditions
(\ref{L2cond0}) and (\ref{L2cond}) in Proposition \ref{prop5} by
$$
q:=\|K_{an}\|_{L^2}^2= \frac13((x^\ast)^3+(1-x^\ast)^3)=\frac13-x^\ast(1-x^\ast)\in [\frac1{12},\frac13].
$$

With this at hand, for an arbitrary sequence of $\alpha=\{\alpha_\omega\}$ satisfying (\ref{L2cond0}) we have
\bea
|f^{(\omega)}_\omega(x_\omega)|^2&\le&\left(\sum_{{\omega'}\subset \omega^c} 
\|f^{(\omega\cup {\omega'})}(x_\omega,t_{\omega'},x^\ast_{(\omega\cup {\omega'})^c})\|_{L_2(X^{\omega'})} q^{|{\omega'}|/2}\right)^2 \qquad (\homega:={\omega}\cup{\omega'})\\
&=&q^{-|\omega|}\left(\sum_{\homega\supset \omega} \|f^{(\homega)}(x_\omega,t_{\homega\backslash\omega},x^\ast_{\homega^c})\|_{L_2(X^{\homega\backslash\omega})} q^{|\homega|/2}\right)^2\\
&\le&q^{-|\omega|}\left(\sum_{\homega\supset \omega} (q^{|\homega|}\alpha_\homega)^{1/2}\|f^{(\homega)}(x_\omega,t_{\homega\backslash\omega},x^\ast_{\homega^c})
\|_{L_2(X^{\homega\backslash\omega})}\cdot \alpha_\homega^{-1/2}\right)^2\\
&\le & C'\alpha_\omega\sum_{\homega \supset \omega} \alpha_\homega^{-1}\|f^{(\homega)}(x_\omega,t_{\homega\backslash\omega},x^\ast_{\homega^c})\|_{L_2(X^{\homega\backslash\omega})}^2,
\eea
where we have replaced $\sum_{\homega\supset \omega} q^{|\homega|}\alpha_\homega$ by $C'q^{|\omega|}\alpha_\omega$ using (\ref{L2cond0}) after applying the Cauchy-Schwarz inequality in the last estimation step.

Now we integrate with respect to $x_\omega\in I^\omega$ and recognize that
$$
\int_{I^\omega} \|f^{(\homega)}(x_\omega,t_{\homega\backslash\omega},x^\ast_{\homega^c})\|_{L_2(X^{\homega\backslash\omega})}^2\,dx_\omega =\|{f'_\homega}^{(\homega)}\|_{L_2}^2.
$$
Thus,
\be\label{Tail}
\|f^{(\omega)}_\omega\|^2_{L_2} \le C'\alpha_\omega \sum_{\homega \supset \omega} \alpha_\homega^{-1}\|{f'_\homega}^{(\homega)}\|_{L_2}^2,\qquad \omega\in
\omega_\gamma,
\ee
and, after substitution and change of summation order, we obtain 
\bea
\|f\|_{\gamma,A}^2 &=& \sum_\omega \frac1{\gamma_\omega}\|f_\omega^{(\omega)}\|_{L_2}^2
\le C'\sum_\omega \frac{\alpha_\omega}{\gamma_\omega} \sum_{\homega\supset \omega} \alpha_\homega^{-1}\|{f'_\homega}^{(\homega)}\|_{L_2}^2\\
&=& C'\sum_\homega \frac1{\gamma_\homega}\|{f'_\homega}^{(\homega)}\|_{L_2}^2 \cdot\left(\frac{\gamma_\homega}{\alpha_\homega}\sum_{\omega\subset \homega} 
\frac{\alpha_\omega}{\gamma_\omega} \right)\le C'C''\|f\|_{\gamma,an}^2,
\eea
where the last step follows from the assumption  (\ref{L2cond}). 

This gives the first half of the proof of the sufficiency of (\ref{L2cond0}). The second half, i.e., the proof that the ANOVA norm dominates the anchored norm, works the same way. Indeed, the formula for ${f'_\omega}^{(\omega)}$ reads now
$$
{f'_\omega}^{(\omega)}(x_\omega)=\sum_{{\omega'}\subset \omega^c} \int_{X^{\omega'}} (P_{(\omega\cup {\omega'})^c}f)^{(\omega\cup {\omega'})}(x_\omega,t_{\omega'})\cdot \kappa_{\omega',A}(x^\ast_{\omega'},t_{\omega'})\,dt_{\omega'},
$$
where $\kappa_{{\omega'},A}(x^\ast_{\omega'},t_{\omega'})=\prod_{k\in {\omega'}} \kappa_{A}(x^\ast,t_k)$. Since $\kappa_{A}(x^\ast,t)=-K_{an}(t)$, the squared $L_2$-norm of
$\kappa_{A}(x^\ast,t)$ has the same value $q$. The remaining estimates are identical. This proves the norm equivalence
\be\label{NE}
C^{-1}\|f\|_{\gamma,an}\le \|f\|_{\gamma,A}\le C\|f\|_{\gamma,an}
\ee
with constant $C:=\sqrt{C'C''}$ if $\gamma$ satisfies the conditions (\ref{L2cond0}) and (\ref{L2cond}).

We next compare the  obtained result with the condition stated in \cite{HS}, where only the case $x^\ast=0$ is detailed.
In the notation adopted in our paper, the condition from \cite{HS} requires the existence of two constants $C',C''$ such that
\be\label{L2cond1}
\sum_{\omega''\supset {\omega'}} 2^{-|\omega''|}\sqrt{\gamma_{\omega''}} \le C' 2^{-|{\omega'}|}\sqrt{\gamma_{\omega'}},\qquad
\sum_{\omega\subset {\omega'}} \frac{1}{\sqrt{\gamma_\omega}}\le \frac{C''}{\sqrt{\gamma_{\omega'}}},\qquad \omega'\in\mP_f(\mathbb{N}).
\ee
Since for $x^\ast=0$ we have $q=1/3$, by choosing $\alpha_\omega=\tilde{q}^{|\omega|}\sqrt{\gamma_{\omega}}$ we see that (\ref{L2cond1}) implies (\ref{L2cond}) for any $1\le \tilde{q} \le 3/2$. 
The proof of Proposition \ref{prop5}
is complete.

\medskip
Even though the conditions given in Proposition \ref{prop5} look more general than (\ref{L2cond1}) it is not easy to construct examples of 
weight sequences $\gamma$ for which (\ref{L2cond0}) and (\ref{L2cond}) hold for some auxiliary sequence $\alpha$ while (\ref{L2cond1}) does not. E.g., for product weights (\ref{prodweights}), one 
easily sees that the two conditions are equivalent, and satisfied if and only if 
\be\label{PWL1}
\sum_{k=1}^\infty \sqrt{\gamma_k} < \infty,
\ee
see \cite{HR,HS}. However, what we have shown is that a direct proof of the norm equivalence (\ref{NE}) yields an at least as good condition as the complex interpolation method exploited in \cite{HS}.

\subsection{Applications to sensitivity analysis}

Here, we consider Sobol indices appearing in connection with high-dimensional model representation techniques and sensitivity analysis for complex systems, where ANOVA and anchored decompositions are considered under various names, see \cite[Chapter 13]{Sm}, \cite{Ra,Ra1}. Given an ANOVA type decomposition $f=\sum_{\omega} f_\omega$ of a multivariate function in a tensor product RKHS, it is customary to measure 
the "influence" of a certain variable group with index set $\omega_0$ together with its interactions with all other variables by the quantity
$$
S_{\omega_0,\mathrm{tot}}:=\sum_{\omega\supset \omega_0} S_{\omega},\qquad S_{\omega}:=\frac{\|f_{\omega}\|^2}{\sum_{\omega'}  \|f_{\omega'}\|^2},
$$
called total Sobol index or total sensitivity index, see \cite[Chapter 15]{Sm}. In the context of weighted tensor product constructions with $H^1(I)$ as coordinate
spaces considered in Section \ref{sec5}, the natural definition of weighted total sensitivity indices would read
\be\label{SobolA}
S_{\omega_0,\mathrm{tot}}^{A,\gamma}:=\sum_{\omega\supset \omega_0} S_{\omega}^{A,\gamma},\qquad S_{\omega}^{A,\gamma}:=\frac{\gamma_{\omega}^{-1} \|f^{(\omega)}_{\omega}\|^2_{L_2}}{\sum_{\omega'} \gamma_{\omega'}^{-1} \|f^{(\omega')}_{\omega'}\|_{L_2}^2},
\ee
for the ANOVA decomposition, and
\be\label{Sobolan}
S_{\omega_0,\mathrm{tot}}^{an,\gamma}:=\sum_{\omega\supset \omega_0} S_{\omega}^{an,\gamma},\qquad S_{\omega}^{an,\gamma}:=\frac{\gamma_{\omega}^{-1} \|{f'}^{(\omega)}_{\omega}\|_{L_2}^2}{\sum_{\omega'} \gamma_{\omega'}^{-1} \|{f'}^{(\omega')}_{\omega'}\|_{L_2}^2},
\ee
for the anchored decomposition, respectively.  The proof of Proposition \ref{prop5} implies the following
\begin{cor}\label{cor2}
Under the conditions (\ref{L2cond0}) and (\ref{L2cond}) on the weight sequence $\gamma$, the total sensitivity indices $S_{\omega_0,\mathrm{tot}}^{A,\gamma}$ and $S_{\omega_0,\mathrm{tot}}^{an,\gamma}$
of any $f$ with finite ANOVA or anchored norm as defined in (\ref{AanNorms}) are comparable:  
$$
C^{-2}\le S_{\omega_0,\mathrm{tot}}^{A,\gamma}/S_{\omega_0,\mathrm{tot}}^{an,\gamma}\le C^2,\qquad \forall\;\omega_0\in\omega_\gamma.
$$
Here, $C$ is the constant in (\ref{NE}).
\end{cor}
The proof is straightforward, just start from (\ref{Tail}) for $\omega\supset \omega_0$, then proceed with the change of summation step
leading to
$$ 
\sum_{\omega\supset \omega_0}\gamma_\omega^{-1}\|f^{(\omega)}_\omega\|_{L_2}^2  =\sum_{\homega\supset\omega_0} \frac1{\gamma_\homega}\|{f'_\homega}^{(\homega)}\|_{L_2}^2 \cdot\left(\alpha_\homega\sqrt{\gamma_\homega}\sum_{\omega_0\subset\omega\subset \homega} 
\frac{\sigma_\omega}{q^{|\omega|}\gamma_\omega} \right)\le C \sum_{\homega\supset\omega_0} \gamma_\homega^{-1}\|{f'_\homega}^{(\homega)}\|_{L_2}^2.
$$
The other direction is analogous. Together with the norm equivalence (\ref{NE}) for $f$, the statement of the corollary follows by forming the corresponding
quotients.  

\smallskip
The message is that a weighted total Sobol index $S_{\omega_0,\mathrm{tot}}^{an,\gamma}$  based on the computationally more feasible anchored decomposition qualitatively carries the same information as the analogous $S_{\omega_0,\mathrm{tot}}^{A,\gamma}$ obtained from the ANOVA type decomposition if a norm equivalence for the underlying spaces $H_{\bfW,\bfa}^\infty$ can be established. This argument fails for the usual unweighted 
Sobol indices based on measuring the terms in the decompositions using the standard $L_2$ norm for which such a norm equivalence cannot hold.

\smallskip
A related question of practical interest is how the approximation error obtained from truncating the anchored decomposition compares to the corresponding
approximation error obtained from truncating the ANOVA decomposition. This problem has been studied in \cite{W,WW1,WW2} for a whole class of weighted tensor product RKHS
and $L_2$ approximation. We consider a particular truncation method called $m$-variate truncation \cite{Ra,Ra1}, similar to (\ref{mVar}), and stick to the concrete setting of Section
\ref{sec5}. Denote by
$$
s^{an}_m(x) := \sum_{|\omega|\le m} f'_\omega(x_\omega), \qquad
s^{A}_m(x) := \sum_{|\omega|\le m} f_\omega(x_\omega),
$$
the $m$-variate approximations obtained from the anchored and ANOVA decompositions of $f$ with finite decomposition norms $\|f\|_{\gamma,an}$ and
$\|f\|_{\gamma,A}$, see (\ref{AanNorms}). 
We assume that
\be\label{GammaL1}
\sum_{\omega} \hat{q}^{|\omega|}\gamma_\omega <\infty
\ee
for some specific value $\hat{q}$ which depends on the decomposition as specified below. The condition (\ref{GammaL1}) does not seem to be too restrictive, for product weights (\ref{prodweights}) it is an automatic consequence of the necessary and sufficient condition (\ref{PWL1})
for the norm equivalence (\ref{NE}). 

By construction of the anchored decomposition we have 
$$
f'_{\omega}(x_{\omega})=((\mathrm{Id}-P')_\omega f(\cdot,x^\ast_{\omega^c}))(x_\omega)=\int_{X^\omega} {f'}^{(\omega)}_\omega(t_\omega)\kappa_{\omega,an}(x_\omega,t_\omega)\, dt_\omega.
$$
Applying the univariate
inequality 
$$
\int_0^1 \left|\int_0^1 g'(t)\kappa_{an}(x,t)\,dt\right|^2 dx\le \int_0^1 g'(t)^2\,dt \int_0^1\int_0^1 \kappa_{an}(x,t)^2\,dtdx =\hat{q}\|g'\|^2_{L_2},
$$
in each coordinate direction appearing in $\omega$, where  $\hat{q}:=\frac12 \max((x^\ast)^2,(1-x^\ast)^2)\in[\frac18,\frac12]$,  we obtain
$$
\|f'_{\omega}\|_{L_2}^2\le  \hat{q}^{|\omega|} \|{f'}^{(\omega)}_{\omega}\|_{L_2}^2.
$$
This allows us to estimate
\bea
\|f-s_m^{an}\|_{L_2}^2&=&\|\sum_{|\omega|>m} f'_\omega\|_{L_2}^2
\le (\sum_{|\omega|>m} \hat{q}^{|\omega|}\gamma_{\omega})
(\sum_{|\omega|>m} \hat{q}^{-|\omega|}\gamma_{\omega}^{-1}\|f'_\omega\|_{L_2}^2)\\
&\le& (\sum_{|\omega|>m} \hat{q}^{|\omega|}\gamma_{\omega})
(\sum_{|\omega|>m} \gamma_{\omega}^{-1}\|{f'}^{(\omega)}_\omega\|_{L_2}^2).
\eea
As a result, we have 
\be\label{anchorrate}
\|f-s_m^{an}\|_{L_2} \le\epsilon_{m,an} \|f\|_{\gamma,an},\qquad \epsilon_{m,an}:=(\sum_{|\omega|>m} \hat{q}^{|\omega|}\gamma_{\omega})^{1/2}\to 0,\quad m\to\infty,
\ee 
due to (\ref{GammaL1}). 

For the terms in the ANOVA decomposition, we similarly have 
$$
f_{\omega}(x_{\omega})=((\mathrm{Id}-P)_\omega P_{\omega^c}f(x_\omega)=\int_{X^\omega} {f}^{(\omega)}_\omega(t_\omega)\kappa_{\omega,A}(x_\omega,t_\omega)\, dt_\omega.
$$
From the associated univariate inequality 
$$
\int_0^1 \left|\int_0^1 g'(t)\kappa_{A}(x,t)\,dt\right|^2 dx\le \int_0^1 g'(t)^2\,dt \int_0^1\int_0^1 \kappa_{A}(x,t)^2\,dtdx=\frac16\|g'\|^2_{L_2},
$$
it follows that
$$
\|f_{\omega}\|_{L_2}^2\le  6^{-|\omega|} \|f^{(\omega)}_{\omega}\|_{L_2}^2.
$$
Repeating the above estimation steps, we obtain
\be\label{ANOVArate}
\|f-s_m^{A}\|_{L_2} \le \epsilon_{m,A} \|f\|_{\gamma,A},\qquad \epsilon_{m,A}:=(\sum_{|\omega|>m} 6^{-|\omega|}\gamma_{\omega})^{1/2}\to 0,\quad m\to\infty,
\ee
for the $m$-variate ANOVA approximation if we assume (\ref{GammaL1}) to hold with $\hat{q}=\frac16$. 

We note that the upper bounds (\ref{anchorrate}) and (\ref{ANOVArate}) are partial cases of \cite[Proposition 1]{W}, in this paper lower bounds are also mentioned. Combining these known results with Proposition \ref{prop5} and using the above notation for weighted Sobol indices, we arrive at
\begin{cor}\label{cor3}
Under the conditions (\ref{L2cond0}) and (\ref{L2cond}) on the weight sequence $\gamma$, the $m$-variate 
anchored approximation $s_m^{an}$ of a function $f$ with finite ANOVA norm $\|f\|_{\gamma,A}<\infty$ satisfies the bound
$$
\|f-s_m^{an}\|_{L_2} \le C\|f\|_{\gamma,A}\left((\sum_{|\omega|>m} \hat{q}^{|\omega|}\gamma_\omega) \sum_{|\omega|>m} S^{A,\gamma}_\omega\right)^{1/2},
$$
where $\hat{q}=\frac12 \max((x^\ast)^2,(1-x^\ast)^2)$ depends on the anchor $x^\ast$, and $C$ is the constant in (\ref{NE}). 
If $1-1/\sqrt{3}<x^\ast<1/\sqrt{3}$, then this error 
bound compares favorably with the bound obtainable when the $m$-variate anchored approximation is replaced by its ANOVA counterpart $s_m^{A}$. 
\end{cor}
The last statement in Corollary 3 is obvious, as for those $x^\ast$ sufficiently close to the center of the interval we have $\hat{q}<\frac16$. The maximal gain occurs if $x^\ast=\frac12$, resulting in an extra factor $(3/4)^{(m+1)/2}$ in the bound for the error of the $m$-variate anchored approximation.

\subsection{Least-squares approximation of maps}\label{sec6}

We briefly discuss the extension of least-squares data fitting or regression problems for functions to the recovery of maps
or operators. Roughly speaking, instead of reconstructing a function $f: X \to \mathbb{R}$ belonging to some hypothesis space from a finite but large
set of noisy samples $(x^i,y^i)$, where $x^i\in X$, and $y^i=f(x^i)+\delta^i$ for $i=1,\ldots,m$, we now want to recover a nonlinear map $F:\, K \to \tK$ from noisy samples $(u^i,\tu^i)\in K\times\tU$ (or their sampled or approximate versions) such that
$$
\tu^i = F(u^i) + {\delta}^i,\qquad i=1,\ldots,n,
$$
with some additive "noise" ${\delta}^i\in \tU$. Here, $K$ and $\tK$ are subsets of certain infinite-dimensional spaces  $U$ and $\tU$, respectively. We are looking for families of reasonable hypothesis spaces our $F$ should belong to. 

Recall first some facts about the nonparametric least-squares regression problem for real-valued functions \cite{CZ,GKKW}. Assuming a generally unknown probability measure $\rho$ on $X\times \mathbb{R}$, 
one tries to "learn" the regression function $f_\rho:\, X\to \mathbb{R}$ defined as the expected value of $y$ given $x\in X$, i.e.,
$$
f_\rho(x)=\int_\mathbb{R} y \,d\rho(y|x),
$$
from a given sample $\bfz=\{x^i,y^i\}_{i=1,\ldots,n}$ which is randomly and independently drawn from $X\times \mathbb{R}$ according to $\rho$. To make this a meaningful problem, it is assumed that $\rho$ is such that $f_\rho$ belongs to a hypothesis Hilbert space $H$ or can be approximated sufficiently well by functions $f:\,X\to \mathbb{R}$ from  $H$. One then encounters the following problem: For  given $\lambda>0$, find $f\in H$ such that the penalized sample error
\be\label{OptL}
      F_\lambda(f):= \frac{1}{n} \sum_{i=1}^n |y^i-f(x^i)|^2 + \lambda \|f\|_H^2 
\ee
is minimized. This approach is also called regularized least-squares regression. Depending on the application, the least-squares data fidelity term 
$\frac1n \sum_{i=1}^n |y^i-f(x^i)|^2$ and the smoothing term $\|f\|_H^2 $ in (\ref{OptL}) can be replaced by more general empirical risk measures and (semi-)norms, respectively.

For certain classes of continuous $f_\rho$ over compact metric spaces $X$, convenient candidates for $H$ are RKHS over $X$ associated with a positive-definite continuous kernel $\kappa:\, X\times X\to \mathbb{R}$. Equivalently, the defining property of a RKHS $H$ is that, for each fixed $x\in X$, the point evaluation $f \to f(x)$ defines a bounded functional on $H$. The kernel $\kappa$ is then given by the formula $\kappa(x,x')=g_x(x')$, where $g_x\in H$ satisfies $f(x)=(f,g_x)_H$ for all $f\in H$ and is unique by the Riesz representation theorem. Besides of the continuity of point evaluations, the main advantage of using RKHS as hypothesis space is that, due to the representer theorem for RKHS, the minimizer of (\ref{OptL}) can be found in form of a finite linear combination
$f(x)=\sum_{i=1}^n c_i \kappa(x,x^i)$ involving the kernel which turns   (\ref{OptL}) into a finite-dimensional quadratic optimization problem.

In order to generalize the sketched kernel-based least-squares method to the reconstruction of maps $F$,  the following approach seems plausible
if the underlying spaces $U$ and $\tU$ possess CONS or Schauder bases, or are equipped with injective parametrizations. Replacing
the given samples $(u^i,\tu^i)$ by the corresponding pairs of sequences of expansion coefficients $(x^i,\tx^i)\in \mX\times \tmX$, we end up with the problem of reconstructing a nonlinear map $\Phi: \mX \to \tmX$ from these derived samples. The coordinate maps $\tx_l=\Phi_l(x)$ of $\Phi$ can be recovered from the data $(x^i,\tx_l^i)\in \mX\times \mathbb{R}$ separately for each $l$ using least-squares regression with suitable spaces $H^\infty_{\bfW,\bfb^l}$ as hypothesis space. Here, the weight sequences $\bfb^l$ can depend on $l$. This approach suffers from the necessity to compute sequences of expansion coefficients $x^i\in \mX$ from the $u^i\in K\subset U$, and $\tx^i \in \tmX$ from the noisy data $\tu^i\in \tU$. Moreover, the implicit assumption that the coordinate maps $\Phi_l$ belong to tensor product Hilbert spaces $H^\infty_{\bfW,\bfb^l}$ for a fixed subspace family $\bfW$ needs justification.

We give some more detail when $\tU$ is a Hilbert space equipped with a CONS $\{\tph_l\}_{l\in\mathbb{N}}$ which provides the parametrization $\tX=\ell^2$. We note that a formal theory of RKHS of Hilbert space valued functions is available \cite{CVT}, the kernel now being a function defined on $\mX\times \mX$ with
values in the space of bounded linear operators from $\tU$ into itself, but we will not make use of it, and stick to our elementary discussion. Then
$$
\tx=\Phi(x)=(\Phi_1(x),\Phi_2(x),\ldots):\qquad \Phi_l(x)=(F(u(x)),\tph_l)_{\tU},\qquad l\in\mathbb{N},
$$
where $u(x)\in U$ is the element represented by $x\in \mX$. Accordingly, $F$ can be recovered from $\Phi$ by summation:
$$
F(u)=F_\Phi(x):=\sum_{l=1}^\infty \Phi_l(x)\tph_l,\qquad u=u(x).
$$ 
Moreover, if the parametrization of $K\subset U$ is given by a Schauder basis 
$\{\phi_k\}_{k\in \mathbb{N}}$ in a Banach space $U$, we also have a linear representation
$$
u(x)=\sum_{k=1}^\infty x_k\phi_k, \qquad x=(x_1,x_2,\ldots)\in\mX,
$$
and can think of $\mX\subset X^\infty=X_1\times X_2\times \ldots$ for some sequence of domains $X_k\subset \mathbb{R}$. The default choice is $X_k=\mathbb{R}$, other choices such as 
$X_k=[a_k,b_k]$  for some sequence of finite numbers $-\infty < a_k<b_k<\infty$ put implicit restrictions on $K$, as $\mX$ must contain the set of all sequences $x$ such that $u(x)\in K$. Then we choose suitable Hilbert spaces $H_k$ of univariate functions $f_k:\,X_k\to \mathbb{R}$ together with their decomposition into
subspaces $W_{j,k}$, $j\in \mJ_k$, and introduce $H^\infty_{\bfW,\bfb^l}$ with a weight sequence $\bfb^l$ as hypothesis space for $\Phi_l$, $l=1,2,\ldots$.
In the spirit of SS-ANOVA, natural candidates for these $H_k=H_{0,k}\oplus H'_k$ would be RHKS split into a small set of finite-dimensional subspaces $W_{j,k}$, $k=0,\ldots,K_k-1$, representing $H_{0,k}$, and an infinite-dimensional RKHS with kernel $H'_k=W_{K_k,k}$ characterized by a kernel $\kappa'_k$, see Section \ref{Ex3}.
In order to arrive at a computationally feasible optimization problem, the weight sequences should overall contain only finitely many 
non-zero entries. Formally, the penalized least-squares regression problem generalizing (\ref{OptL}) reads then as follows: 
Find $\Phi=(\Phi_1(x),\Phi_2(x),\ldots)$ with $\Phi_l\in H^\infty_{\bfW,\bfb^l}$,
such that
\be\label{FA}
\frac1n \sum_{i=1}^n \|\tu^i-F(u^i)\|_{\tU}^2 + \lambda \sum_{l=1}^\infty \|\Phi_l\|_{\bfW,\bfb^l}^2
\quad \longmapsto \quad \min.
\ee
Explicitly, we have
\be\label{Data}
\|\tu^i-F(u^i)\|_{\tU}^2=\sum_{l=1}^\infty |(\tu^i,\tph_l)_{\tU} -\Phi_l(x^i)|^2 = \sum_{l=1}^\infty |(\tu^i,\tph_l)_{\tU} -\sum_{\bfj\in\mJ} \ww_{\bfj}^l(x^i)|^2, 
\ee
where $\Phi_l(x)=\sum_{\bfj\in\mJ} \ww_{\bfj}^l(x)$, $\ww_{\bfj}^l\in W_{\bfj}$, is the unique decomposition with  of $\Phi_l\in H^\infty_{\bfW,\bfb^l}$
with respect to $\bfW$. Similarly,
\be\label{Smooth}
\sum_{l=1}^\infty \|\Phi_l\|_{\bfW,\bfb^l}^2 = \sum_{l=1}^\infty \sum_{\bfj\in\mJ} b^l_{\bfj}\|\ww_{\bfj}^l\|_{\omega_{\bfj}}^2.
\ee
If only finitely many $b_{\bfj}^l$ are non-zero, then $\Phi_l=0$ for all but finitely many $l$, too, and the summations in both (\ref{Data}) and (\ref{Smooth})
become finite, without changing the minimizer in (\ref{FA}). As is the case for kernel-based regression, and in particular for SS-ANOVA \cite{Gu,Wa}, with this assumption it automatically follows that (\ref{FA}) can be solved as a finite-dimensional quadratic optimization problem using linear combinations
of tensor products of basis functions from the active $W_{j,k}$ with $j<K_k$, and kernels $\kappa'_k$ if some $H'_k=W_{K_k,k}$ is active. Here, active means that the corresponding $b_{\bfj}^l$ does not vanish for at least one $l$. 

The challenge will be to determine the supports $\omega_{\bfb^l}$ and non-zero values $b^l_{\bfj}$ that guarantee interpretable results and error rates
depending on the nature of $K\subset U$ and $\tK\subset \tU$, similar to what is done in statistical learning theory.
If $U$ and $\tU$ are function spaces, an additional issue comes into play: In practice, we do not have access to the samples $(u^i,\tu^i)$ but only to finitely many point values or linear functionals evaluated from $u^i$ and $\tu^i$ (we call these the sampled versions of the samples). How to include this secondary
level of sampling into the above functional fitting scheme needs to be considered for a given application. We believe that a more detailed study of least-squares approximation of maps between infinitely parameterized domains from sampled information is warranted but leave it for future research.


\section*{Acknowledgements} 
M. Griebel acknowledges the support of the Deutsche Forschungsgemeinschaft (DFG) through the Sonderforschungsbereich 1060: The Mathematics of Emergent Effects.
This paper was written while P. Oswald held a Bonn Research Chair sponsored by the
Hausdorff Center for Mathematics at the University of Bonn funded by the Deutsche Forschungsgemeinschaft. He is grateful for this support.

\end{document}